\documentclass[aos,nameyear,dvips]{arximspdf}
\usepackage{dcolumn}
\usepackage{graphicx}

\doi{10.1214/10-AOS812}
\volume{38}
\issue{5}
\pubyear{2010}
\firstpage{2958}
\lastpage{2997}

\makeatletter

\newcolumntype{d}[1]{D{.}{.}{#1}}

\newcommand{\tr}{\operatorname{tr}}
\newcommand{\cum}{\operatorname{cum}}
\newcommand{\ve}{\varepsilon}
\newcommand{\imag}{i}
\newcommand{\arctanh}{\operatorname{arctanh}}
\newcommand{\wihat}{\widehat}
\newcommand{\esp}{E}
\newcommand{\cov}{\operatorname{cov}}
\newcommand{\var}{\operatorname{Var}}

\newproclaim{exmp}{Example}[section]
\newtheorem{theorem}{Theorem}[section]
\newproclaim{rem}{Remark}[section]
\newtheorem{proposition}{Proposition}
\newtheorem{lemma}{Lemma}

\makeatother

\begin{document}
\begin{frontmatter}

\title{An efficient estimator for locally stationary Gaussian
long-memory processes\thanksref{T1}}
\runtitle{Locally stationary long-memory processes}

\thankstext{T1}{Supported in part by Fondecyt Grant 1085239.}

\begin{aug}
\author[A]{\fnms{Wilfredo} \snm{Palma}\corref{}\ead[label=e1]{wilfredo@mat.puc.cl}} and
\author[A]{\fnms{Ricardo} \snm{Olea}}
\runauthor{W. Palma and R. Olea}
\affiliation{Pontificia Universidad Cat\'olica de Chile}
\address[A]{Department of Statistics\\
Pontificia Universidad Cat\'olica de Chile\\
Vicu\~{n}a Mackenna 4860\\
Macul, Santiago\\
Chile\\
\printead{e1}} 
\end{aug}

\received{\smonth{4} \syear{2009}}
\revised{\smonth{1} \syear{2010}}

%
\begin{abstract}
This paper addresses the estimation of locally stationary long-range
dependent processes, a methodology that allows the statistical analysis
of time series data exhibiting both nonstationarity and strong
dependency. A time-varying parametric formulation of these models is
introduced and a Whittle likelihood technique is proposed for
estimating the parameters involved. Large sample properties of these
Whittle estimates such as consistency, normality and efficiency are
established in this work. Furthermore, the finite sample behavior of
the estimators is investigated through Monte Carlo experiments. As a
result from these simulations, we show that the estimates behave well
even for relatively small sample sizes.
\end{abstract}

%

\end{frontmatter}

\section{Introduction}\label{sec1}
Even though stationarity is a very attractive theoretical assumption,
in practice most time series data fail to meet this condition. As a
consequence, several approaches to deal with nonstationarity have been
proposed in the literature. Among these methodologies, differentiation
and trend removal are popular choices. Other approaches include, for
instance, the evolutionary spectral techniques first discussed by
\citet{Prie65}. In a similar spirit, during the last decades a
number of new time-varying dependence models have been proposed. One of
these methodologies, the so-called locally stationary processes
developed by Dahlhaus (\citeyear{Dahl96b}, \citeyear{Dahl97}), has been widely discussed in
the recent time series literature, see, for example,
\citet{Dahl00}, \citet{von00}, \citet{Jens00},
\citet{Guo03}, \citet{gent04}, \citet{orbe05},
Dahlhaus and Polonik (\citeyear{Dahl06}, \citeyear{Dahl09}), \citet{chand06}, \citet{fryz06}
and \citet{Bera09}, among others. This approach allows the
stochastic process to be nonstationary, but assuming that the time
variation of the model is sufficiently smooth so that it can be locally
approximated by stationary processes.

On the other hand, during the last decades, long-range dependent data
have arisen in disciplines as diverse as meteorology, hydrology,
economics, etc., see, for example, the recent surveys by
\citet{Douk03} and \citet{Palm07}. As a consequence,
statistical methods for modeling that type of data are of great
interest to scientists and practitioners from many fields. At the same
time, many of these long-memory data also display nonstationary
behavior, see, for instance, \citet{Gran96}, \citet{Jens00}
and \citet{Bera09}. Nevertheless, most of the currently available
methods for dealing with long-range dependence are incapable of
modeling time series with these features. In particular, much of the
theory of locally stationary processes applies only to time series with
short memory, such as time-varying autoregressive moving average (ARMA)
processes and not to time series exhibiting both nonstationarity and
strong dependence. In order to treat that type of data, this paper
addresses a class of strongly dependent locally stationary processes.
In particular, these models include a Hurst parameter which evolves
over time. Following \citet{Dahl97}, we propose a Whittle maximum
likelihood estimation technique for fitting Gaussian long-memory
locally stationary models. This is an extension of the spectrum-based
likelihood estimator introduced by \citet{Whit53}. A great
advantage of this estimation procedure is its computational efficiency,
since it only requires the calculation of the periodogram by means of
the fast Fourier transform. Additionally, we prove in this article that
the proposed Whittle estimator is asymptotically consistent, normally
distributed and efficient. Thus, this paper provides a framework for
modeling and making statistical inferences about several types of
nonstationarities that may be difficult to handle with other
techniques. For instance, changes in the variance of a time series
could be spotted by simple inspection of the data. However, variations
on the dependence structure of the data are far more difficult to
uncover and model.

The remainder of this paper is structured as follows. Section
\ref{methodology} discusses a class of long-memory locally stationary
processes and proposes a quasi maximum likelihood estimator based on an
extended version of the Whittle spectrum-based methodology. The
consistency, asymptotic normality and efficiency of these quasi maximum
likelihood estimators are established. Applications of the asymptotic
results to some specific locally stationary processes are also
presented in this section. Proofs of the theorems are provided in
Section \ref{proofs}. Note that the techniques employed by
\citet{Dahl97} to show the asymptotic properties of the Whittle
estimates are no longer valid for the class of long-memory locally
stationary processes discussed in this paper. This difficulty is due to
the fact that these processes have an unbounded time-varying spectral
density at zero frequency. Consequently, several technical results must
be introduced and proved. Section \ref{simulations} reports the results
from several Monte Carlo experiments which allow to gain some insight
into the finite sample behavior of the Whittle estimates. Conclusions
are presented in Section \ref{conclusion} while auxiliary lemmas are
provided in a technical appendix. Additional examples and simulations
along with a comparison of the Whittle estimator with a kernel maximum
likelihood estimation approach and two real-life applications of the
proposed methodology can be found in \citet{Palm10a}. The
bandwidth selection problem for the locally stationary Whittle
estimator is also discussed in that paper, from an empirical
perspective.

\section{Definitions and main results} \label{methodology}
\subsection{Long-memory locally stationary processes} A class of
Gaussian locally stationary process with transfer function $A^0$
can be defined by the spectral representation
%
\begin{equation}\label{local-stat}
Y_{t,T} = \int_{-\pi}^{\pi}A^{0}_{t,T}(\lambda) e^{i\lambda t}\,
dB(\lambda),
\end{equation}
for $t=1,\ldots,T$, where $B(\lambda)$ is a Brownian motion
on $[-\pi,\pi]$ and there is a positive
constant $K$ and a $2\pi$-periodic function $A\dvtx (0,1]\times\mathbb
{R}\to
\mathbb{C}$ with $A(u,-\lambda)=\overline{A(u,\lambda)}$ such that
%
\begin{equation} \label{ineq}
\sup_{t,\lambda} \biggl|A^{0}_{t,T}(\lambda)-A \biggl( \frac{t}{T},\lambda\biggr)
\biggr|\leq\frac{K}{T},
\end{equation}
for all $T$. The transfer function $A^{0}_{t,T}(\lambda)$ of this
class of
nontstationary processes changes smoothly over time so that they can be locally
approximated by stationary processes. An example of this class of
locally stationary
processes is given by the infinite moving average expansion
%
\begin{equation} \label{wold-general}
Y_{t,T}=\sigma\biggl( \frac{t}{T} \biggr)\sum_{j=0}^\infty
\psi_j \biggl( \frac{t}{T} \biggr)\varepsilon_{t-j},
\end{equation}
where $\{\varepsilon_t\}$ is a zero-mean and unit variance Gaussian
white noise
and $\{\psi_j(u)\}$ are coefficients satisfying $\sum_{j=0}^\infty
\psi_j (u )^2<\infty$ for all $u \in[0,1]$.
In this case, the transfer function of process (\ref{wold-general}) is
given by
$A^{0}_{t,T}(\lambda)=\break\sigma( \frac{t}{T} )\sum_{j=0}^\infty
\psi_j ( \frac{t}{T} )e^{-i\lambda j}=A ( \frac{t}{T},\lambda)$, so that
condition (\ref{ineq}) is satisfied. The model defined by (\ref
{wold-general}) generalizes
the Wold expansion for a linear stationary process allowing the
coefficients of the
infinite moving average expansion vary smoothly over time. 
A particular case of (\ref{wold-general}) is the generalized version
of the
fractional noise process described by the discrete-time equation
%
\begin{equation} \label{wold-local}
Y_{t,T}=\sigma\biggl( \frac{t}{T} \biggr)(1-B)^{-d ( {t/T} )}\varepsilon_t
=\sigma\biggl( \frac{t}{T} \biggr) \sum_{j=0}^{\infty}
\eta_j \biggl( \frac{t}{T} \biggr) \varepsilon_{t-j},
\end{equation}
for $t=1,\ldots,T$, where $\{\varepsilon_t\}$ is a Gaussian white
noise sequence with
zero mean and unit variance and the infinite moving average
coefficients $\{\eta_j(u)\}$
are given by
%
\begin{equation} \label{ls-fn}
\eta_j (u ) =
\frac{\Gamma[j+d (u )]}{\Gamma(j+1)\Gamma[d (u )]},
\end{equation}
where $\Gamma(\cdot)$ is the Gamma function and $d(\cdot)$ is a
smoothly time-varying
long-memory coefficient. For simplicity, the locally stationary
fractional noise process
(\ref{wold-local}) will be denoted as LSFN.

A natural extension of the LSFN model is the locally stationary autoregressive
fractionally integrated moving average (LSARFIMA) process defined by
the equation
%
\begin{equation} \label{ext-ARFIMA}
\Phi(t/T,B)Y_{t,T}=\sigma(t/T) \Theta(t/T,B) (1-B)^{-d(t/T)} \ve_{t},
\end{equation}
%
for $t=1,\ldots,T$, where for $u\in[0,1]$, $\Phi(u,B)=1+\phi
_1(u)B+\cdots+
\phi_P(u)B^P$ is an autoregressive polynomial, $\Theta(u,B)=1+\theta
_1(u)B+\cdots+
\theta_Q(u)B^Q$ is a moving average polynomial, $d(u)$ is a
long-memory parameter,
$\sigma(u)$ is a noise scale factor and $\{\ve_t\}$ is a Gaussian
white noise sequence
with zero mean and unit variance. This class of models extends the
well-known ARFIMA
process, which is obtained when the components $\Phi(u,B)$, $\Theta
(u,B)$, $d(u)$ and
$\sigma(u)$ appearing in (\ref{ext-ARFIMA}) do not depend on $u$.
Note that by
Theorem 4.3 of \citet{Dahl96b}, under some regularity conditions
on the polynomial
$\Phi(u,B)$, the model defined by (\ref{ext-ARFIMA}) satisfies (\ref
{local-stat}) and~(\ref{ineq}), see \citet{Jens00} for details.

\subsection{Estimation} \label{section-asymp} 
Let $\theta\in\Theta$ be a parameter vector specifying model (\ref
{local-stat}) where
the parameter space $\Theta$ is a subset of a finite-dimensional
Euclidean space. Given a sample $\{Y_{1,T},\ldots,Y_{T,T}\}$ of the
process (\ref{local-stat}) we can
estimate $\theta$ by minimizing the
Whittle log-likelihood function
%
\begin{equation} \label{whittle-like}
\mathcal{L}_T(\theta)=\frac{1}{4\pi}\frac1M \int_{-\pi}^{\pi
}\sum_{j=1}^M \biggl\{ \log
f_{\theta}(u_j,\lambda)+
\frac{I_N(u_j,\lambda)}{f_{\theta}(u_j,\lambda)} \biggr\}\,d\lambda,
\end{equation}
where $f_{\theta}(u,\lambda)=|A_{\theta}(u,\lambda)|^2$ is the
time-varying spectral density of the limiting process specified by the
parameter $\theta$, $I_N(u,\lambda)=\frac{|D_N(u,\lambda)|^2}{2\pi
H_{2,N}(0)}$ is a tapered
periodogram with
\[
D_N(u,\lambda)=\sum_{s=0}^{N-1}h \biggl(\frac{s}{N} \biggr)Y_{[uT]-N/2+s+1,T}
e^{-i\lambda s},\qquad H_{k,N}=\sum_{s=0}^{N-1}h \biggl(\frac{s}{N} \biggr)^k
e^{-i\lambda s},
\]
$T=S(M-1)+N$, $u_j=t_j/T$, $t_j=S(j-1)+N/2$, $j=1,\ldots,M$
and $h(\cdot)$ is a data taper. The intuition behind this extended
version of the Whittle estimation procedure (\ref{whittle-like}) is as
follows: the sample
$\{Y_{1,T},\ldots,Y_{T,T}\}$ is subdivided into $M$ blocks of length
$N$ each shifting $S$ places from block to block. For instance, if we
split a time series of $T=652$ observations into $M=100$ blocks of
length $N=256$ each,
shifting $S=4$ positions forward each time we get the blocks
$(Y_{1,652},Y_{2,652},\ldots,Y_{256,652}), 
\ldots, (Y_{397,652},Y_{398,652},\ldots,Y_{652,652})$. Then, the
spectrum is locally
estimated by means of the data tapered periodogram on each one of these
$M=100$ blocks
and then averaged to form (\ref{whittle-like}). Finally, the Whittle
estimator of the
parameter vector $\theta$ is given by
%
\begin{equation} \label{theta}
\widehat{\theta}_T=\arg\min\mathcal{L}_T(\theta),
\end{equation}
where the minimization is over a parameter space $\Theta$. The
analysis of the
asymptotic properties of the Whittle locally stationary estimates
(\ref{theta}) is
discussed in detail next. Before stating these results, we introduce a
set of the regularity conditions.

\subsection{Assumptions}
The first assumption below is concerned with the time-varying
spectral density of the process. The second assumption is related to
the data tapering
function and the third assumption is concerned with the block sampling
scheme. It is
assumed that the parameter space $\Theta$ is compact. In what follows,
$K$ is always a
positive constant that could be different from line to line.
\begin{longlist}[A3.]
\item[A1.] \hypertarget{A1}
The time-varying spectral density of the limiting process (\ref
{local-stat}) is
strictly positive and satisfies
\[
f_{\theta}(u,\lambda) \sim C_f(\theta, u) |\lambda|^{-2d_{\theta}(u)},
\]
as $|\lambda|\to0$, where $C_f(\theta, u)>0$, $0<\inf_{\theta
,u}d_{\theta}(u)$, $\sup_{\theta,u}d_{\theta}(u) < \frac12$ and
$d_{\theta}(u)$ has bounded first derivative with respect to $u$.
There is an integrable function $g(\lambda)$ such that $|\nabla
_{\theta} \log
f_{\theta}(u,\lambda)|\leq g(\lambda)$ for all $\theta\in\Theta$,
$u\in[0,1]$ and
$\lambda\in[-\pi,\pi]$. The function $A(u,\lambda)$ is
twice differentiable with respect to $u$ and satisfies
\[
\int_{-\pi}^{\pi} A(u,\lambda) A(v,-\lambda)\exp(\imag k \lambda
) \,d\lambda\sim C(\theta,
u,v) k^{d_{\theta}(u)+d_{\theta}(v)-1},
\]
as $k\to\infty$, where $|C(\theta,u,v)|\leq K$ for $u,v\in[0,1]$
and $\theta\in\Theta$.
The function $f_{\theta}(u,\lambda)^{-1}$ is twice differentiable
with respect to
$\theta$, $u$ and $\lambda$.
\item[A2.]  \hypertarget{A3}
The data taper $h(u)$ is a positive, bounded function for $u\in[0,1]$
and symmetric
around $ \frac12$ with a bounded derivative.
\item[A3.]  \hypertarget{A4}
The sample size\vspace*{2pt} $T$ and the subdivisions integers $N$, $S$ and $M$
tend to infinity satisfying $S/N\to0$, $\sqrt{T}\log^2 N/N\to0$, $\sqrt{T}/M\to0$
and $N^3\log^2
N/\break T^2\to0$.
\end{longlist}
\begin{exmp}
As an illustration of the assumptions described above, consider the
extension of the
usual fractional noise process with time-varying Hurst parameter,
described by
(\ref{wold-local}) and (\ref{ls-fn}). The spectral density of this LSFN
process is given
by
\[
f_{\theta}(u,\lambda) = \frac{\sigma^2}{2\pi} \biggl(2 \sin\frac
{\lambda}{2} \biggr)^{-2
d_{\theta}(u)}.
\]
Note that this function is integrable over $\lambda\in[-\pi,\pi]$
for every
$u\in[0,1]$ as long as $d_{\theta}(u)<\frac12$ for all $u \in[0,1]$
and $\theta\in
\Theta$. Furthermore, we have that $f_{\theta}(u,\lambda)\sim\frac
{\sigma^2}{2\pi}
|\lambda|^{-2 d_{\theta}(u)}$, as $\lambda\to0$. By assuming that
$|\nabla_{\theta}
d_{\theta}(u)|\leq K$, we have that $|\nabla_{\theta} \log
f_{\theta}(u,\lambda)|=|\nabla_{\theta}
d_{\theta}(u)||{\log}(2\sin\frac{\lambda}{2} )^2|\leq K |{\log}
|\lambda||$, which
is an integrable function in $\lambda\in[-\pi,\pi]$. In addition,
from (\ref{ls-fn}) the
function $A(u,\lambda)$ of this process satisfies
\[
\int_{-\pi}^{\pi} A(u,\lambda) A(v,-\lambda)\exp(\imag k \lambda
) \,d\lambda=
\frac{\Gamma[1-d_{\theta}(u)-d_{\theta}(v) ]
\Gamma[k+d_{\theta}(u) ]}{\Gamma[1-d_{\theta}(u) ]
\Gamma[d_{\theta}(u) ]\Gamma[k+1-d_{\theta}(v) ]},
\]
for $k\geq0$. Thus, by Stirling's approximation, we get
\[
\int_{-\pi}^{\pi} A(u,\lambda) A(v,-\lambda)\exp(\imag k \lambda
) \,d\lambda
\sim\frac{\Gamma[1-d_{\theta}(u)-d_{\theta}(v) ]}{\Gamma
[1-d_{\theta}(u) ]
\Gamma[d_{\theta}(u) ]}k^{d_{\theta}(u)+d_{\theta}(v)-1},
\]
for $k\to\infty$. Besides, a simple calculation shows that $f_{\theta
}(u,\lambda)^{-1}$
is twice differentiable with respect to $u$ and $\lambda$ as long as
$d_{\theta}(u)$ is
twice differentiable with respect to $u$. Thus, under these conditions
the time-varying
spectral density $f_{\theta}(u,\lambda)$ satisfies assumption \hyperlink{A1}{A1}. On the other
hand, an example of data taper that satisfies assumption \hyperlink{A3}{A2} is
the cosine bell
function
%
\begin{equation} \label{bell}
h(x)= \tfrac12[1-\cos(2\pi x)].
\end{equation}
Note that if $S=\mathcal{O}(N^{a})$ and $M=\mathcal{O}(N^{b})$ then
$T=\mathcal{O}(N^{a+b})$ for $a+b\geq1$. Thus, by choosing exponents
$a$ and $b$ such
that $(a,b) \in\mathcal{C}=\{a<1, \frac32<a+b<2, a<b\}$, assumption
\hyperlink{A4}{A3} is fulfilled. Observe that the $\mathcal{C}$ is a nonempty set.
\end{exmp}

\subsection{Main results} \label{main-results}
Some fundamental large sample properties of the Whittle
quasi-likelihood estimators
(\ref{theta}), including consistency, asymptotic normality and
efficiency are established next. In addition, we establish an
asymptotic result about the estimation of the time-varying long-memory
parameter for a class of locally stationary processes. The proofs of
these four results are provided in Section \ref{proofs}.
\begin{theorem}[(Consistency)] \label{CON} Let $\theta_0$ be the true value
of the parameter $\theta$.
Under assumptions \textup{\hyperlink{A1}{A1}--\hyperlink{A4}{A3}}, the estimator $\wihat{\theta
}_T$ satisfies
$\wihat{\theta}_T \to\theta_0$, in probability, as $T\to\infty$.
\end{theorem}
\begin{theorem}[(Normality)] \label{CLT} Let $\theta_0$ be the true value
of the parameter $\theta$.
If assumptions \textup{\hyperlink{A1}{A1}--\hyperlink{A4}{A3}} hold, then the Whittle estimator
$\wihat{\theta}_T$
satisfies a central limit theorem
\[
\sqrt{T}( \wihat{\theta}_T-\theta_0) \to N [0,\Gamma(\theta
_0)^{-1} ],
\]
in distribution, as $T\to\infty$, where
%
\begin{equation} \label{Gamma}
\Gamma(\theta)= \frac{1}{4\pi} \int_{0}^{1}\int_{-\pi}^{\pi}
[\nabla\log
f_{\theta}(u,\lambda)] [\nabla\log
f_{\theta}(u,\lambda)]' \,d\lambda \,du.
\end{equation}
\end{theorem}
\begin{theorem}[(Efficiency)] \label{EFF} Assuming that conditions
\textup{\hyperlink{A1}{A1}--\hyperlink{A4}{A3}} hold, the Whittle estimator $\wihat{\theta}_T$
is asymptotically Fisher efficient.
\end{theorem}
\begin{rem}
Recall that for a stationary fractional noise process FN($d$), the
asymptotic variance of the
maximum likelihood estimate of the long-memory parameter, $\widehat
{d}$, satisfies
\[
\lim_{T\to\infty} T\var(\widehat{d})=\frac{6}{\pi^2}.
\]
On the other hand, suppose that we consider a LSFN process where the long-memory
parameter varies according to, for example, $d(u)=\alpha_0+\alpha_1
u$. Thus, in order
to estimate $d(u)$, the parameters $\alpha_0$ and $\alpha_1$ must be
estimated. Let $\widehat{\alpha}_0$ and $\widehat{\alpha}_1$ be
their Whittle estimators, respectively, so that $\widehat
{d}(u)=\widehat{\alpha}_0+\widehat{\alpha}_1 u$. According to
Theorem \ref{CLT}, the asymptotic variance of this estimate of $d(u)$ satisfies
\[
\lim_{T\to\infty} T\var[\widehat{d}(u)]=\frac{24}{\pi^2}(1-3u+3u^2),
\]
and then integrating over $u$ we get
\[
\lim_{T\to\infty} T\int_0^1\var[\widehat{d}(u)]\,du=\frac{12}{\pi^2}.
\]
Since two parameters are being estimated, on the average, the
asymptotic variance of the
estimate $\widehat{d}(u)$ is twice the asymptotic variance of
$\widehat{d}$ from a stationary FN process. This result can be
generalized to the case where three or more coefficients are estimated
and to more complex trends, as established on the following theorem.
\end{rem}
\begin{theorem} \label{average-d}
Consider a LSFN process (\ref{wold-local}) with time-varying
long-memory parameter
$d_{\beta}(u) = \sum_{j = 1}^{p} \beta_j g_j(u)$, where $\{g_j(u)\}$
are basis functions
as defined in (\ref{theta-g}) below. Let $\wihat{d}(u) = \sum_{j =
1}^{p} \wihat{\beta}_j
g_j(u)$ be the estimator of $d_{\beta}(u)$ for $u\in[0,1]$. Then
under assumptions
\textup{\hyperlink{A1}{A1}--\hyperlink{A4}{A3}} we have that
%
\begin{equation} \label{ave-limit}
\lim_{T\to\infty} T \int_{0}^{1} \var[ \wihat{d}(u) ] \,du =
\frac{6 p}{\pi^2}.
\end{equation}
\end{theorem}
\begin{rem}
Note that according to Theorem \ref{average-d} the limiting average of
the variances of
$d(u)$ given by (\ref{ave-limit}) does not depend on the basis
functions $g_j(\cdot)$
for $j=1,\ldots,p$.
\end{rem}

\subsection{Illustrations} \label{examples} As an illustration of the
asymptotic results discussed above, consider the class of LSARFIMA
models defined by
(\ref{ext-ARFIMA}). The evolution of these models
can be specified in terms of a general class of functions. For example,
let $\{g_j(u)\}$, $j=1,2,\ldots,$ be a basis for a space of smoothly varying
functions and let $d_{\theta}(u)$ be the time-varying long-memory
parameter in model (\ref{ext-ARFIMA}). Then we could write $d_{\theta
}(u)$ in terms of the basis $\{g_j(u)\}$ as follows:
%
\begin{equation} \label{theta-g}
\ell[d_{\theta}(u)]=\sum_{j=0}^k \alpha_j g_j(u),
\end{equation}
for unknown values of $k$ and $\theta=(\alpha_0,\alpha_1,\ldots
,\alpha_k)'$, where $\ell(\cdot)$ is a known link function. In this
situation, estimating $\theta$ involves determining $k$ and estimating
the coefficients $\alpha_0,\ldots,\alpha_k$. Important examples of
this approach are the classes of polynomials generated by the basis $\{
g_j(u)=u^j\}$,
Fourier expansions generated by the basis $\{g_j(u)=e^{\imag u j}\}$
and wavelets
generated by, for instance, the Haar or Daubechies systems. Extensions
of these cases
can also be considered. For example, the basis functions could also
include parameters
as in the case $\{g_j(u)=e^{\imag u \beta_j }\}$, where $\{\beta_j\}$
are unknown
values.

In order to illustrate the application of the theoretical results
established in
Section \ref{main-results}, we discuss next a number of
combinations of polynomial and
harmonic evolutions of the long-memory parameter, the noise variance, the
autoregressive and moving average components of the LSARFIMA process
(\ref{ext-ARFIMA}).
Additional examples are provided in Section 2 of \citet{Palm10a}.
\begin{exmp} \label{examp-general}
Consider first the case $P=Q=0$ in model (\ref{ext-ARFIMA}) where
$d(u)$ and $\sigma(u)$ are specified by
\[
\ell_1 [ d(u) ] = \sum_{j=0}^p \alpha_j g_j(u),\qquad \ell_2 [ \sigma(u)
] = \sum_{j=0}^q \beta_j h_j(u),
\]
for $u\in[0,1]$, where $\ell_1(\cdot)$ and $\ell_2(\cdot)$ are
differentiable link
functions, $g_j(\cdot)$ and $h_j(\cdot)$ are basis functions. The
parameter vector in
this case is $\theta=(\alpha_0,\ldots,\alpha_p$, $\beta_0,\ldots
,\beta_q)'$ and the matrix $\Gamma$
can be written as
%
\begin{equation} \label{Gamma-a-b}
\Gamma= \pmatrix{
\Gamma_{\alpha} & 0 \cr
0 & \Gamma_{\beta} },
\end{equation}
where
\begin{eqnarray*}
\Gamma_{\alpha} &=& \frac{\pi^2}{6 } \biggl[
\int_0^1 \frac{g_i(u) g_j(u) \,du }{[ \ell_1^{\prime}(d(u)) ]^2}
\biggr]_{i,j=0,\ldots,p} ,\\
\Gamma_{\beta} &=& 2 \biggl[
\int_0^1 \frac{h_i(u) h_j(u) \,du}{ [ \sigma(u) \ell_2^{\prime}(\sigma
(u)) ]^2}
\biggr]_{i,j=0,\ldots,q}.
\end{eqnarray*}
\end{exmp}
\begin{exmp} \label{examp-poly}
As a particular case of the parameter specification of the previous
example, consider
the case $P=Q=0$ in model (\ref{ext-ARFIMA}) where
$d(u)$ and $\sigma(u)$ are both specified by polynomials,
\[
d(u) = \alpha_0 + \alpha_1 u + \cdots+ \alpha_p u^p,\qquad
\sigma(u) = \beta_0 + \beta_1 u + \cdots+ \beta_q u^q,
\]
for $u\in[0,1]$. Similar to Example \ref{examp-general}, in this
case the parameter
vector is $\theta=(\alpha_0,\ldots,\alpha_p, \beta_0,\ldots,\beta
_q)'$, $\ell_1(u)=\ell_2(u)=u$ and the matrix $\Gamma$
given by (\ref{Gamma}) can be written as in (\ref{Gamma-a-b}) with
\begin{eqnarray*}
\Gamma_{\alpha} &=& \biggl[\frac{\pi^2}{6 (i+j+1)}
\biggr]_{i,j=0,\ldots,p},\\
\Gamma_{\beta}&=&2 \biggl[
\int_0^1 \frac{u^{i+j} \,du}{(\beta_0 + \beta_1 u + \cdots+ \beta_q
u^q)^2} \biggr]_{i,j=0,\ldots,q}.
\end{eqnarray*}
The above integrals can be evaluated by standard calculus procedures;
see, for example,
Gradshteyn and Ryzhik [(\citeyear{Grad00}), page 64] or by numerical integration.
\end{exmp}
\begin{exmp} \label{examp-exp}
Considering now a similar setup as Example \ref{examp-poly} with
$p=q=1$, but with link function $\ell(\cdot)=\log(\cdot)$ such that
\[
\log[ d(u) ] = \alpha_0 + \alpha_1 u,\qquad
\log[ \sigma(u) ] = \beta_0 + \beta_1 u,
\]
for $u\in[0,1]$. Then $\Gamma$ can be written as (\ref{Gamma-a-b}) with
\begin{eqnarray*}
\Gamma_{\alpha} &=& \frac{\pi^2}{6}\frac{e^{2\alpha_0}}{4\alpha
_1^3} \left[ \matrix{
2\alpha_1^{2}(e^{2\alpha_1}-1) & \alpha_1\bigl((2\alpha_1-1)e^{2\alpha_1}+1\bigr)
\vspace*{2pt}\cr
\alpha_1\bigl((2\alpha_1-1)e^{2\alpha_1}+1\bigr) &
(2\alpha_1^2-2\alpha_1+1)e^{2\alpha_1}+1}
\right],
\\
\Gamma_{\beta} &=& \left[ \matrix{
2 & 1 \cr
1 & 2/3}
\right].
\end{eqnarray*}
\end{exmp}
\begin{exmp} Following with the assumption $P=Q=0$ in model (\ref
{ext-ARFIMA}), consider
that $d(u)$ and $\sigma(u)$ are defined by the harmonic expansions
\begin{eqnarray*}
d(u) & = & \alpha_0 + \alpha_1 \cos(\lambda_1 u) + \cdots+ \alpha_p
\cos(\lambda_p u),\\
\sigma(u) & = & \beta_0 + \beta_1 \cos(\omega_1 u) + \cdots+ \beta
_q \cos(\omega_q u),
\end{eqnarray*}
for $u\in[0,1]$, where $\lambda_0=0$, $\lambda_i^2\neq\lambda_j^2$
for all
$i,j=0,\ldots,p$, $i\neq j$, $\omega_0=0$ and $\omega_i^2\neq\omega
_j^2$ for all
$i,j=0,\ldots,q$, $i\neq j$. For simplicity, the values of the frequencies
$\{\lambda_j\}$ and $\{\omega_j\}$ are assumed to be known. As in Example
\ref{examp-poly}, in this case the parameter vector is $\theta
=(\alpha_0,\ldots,\alpha_p,
\beta_0,\ldots,\beta_q)'$ and the matrix $\Gamma$ appearing in
(\ref{Gamma}) can be
written as in (\ref{Gamma-a-b}) with
\[
\Gamma_{\alpha}=\frac{\pi^2}{12}
\biggl[\frac{\sin(\lambda_i-\lambda_j)}{\lambda_i-\lambda_j}
+ \frac{\sin(\lambda_i+\lambda_j)}{\lambda_i+\lambda_j}
\biggr]_{i,j=0,\ldots,p}
\]
and
\[
\Gamma_{\beta}=\frac{\pi^2}{12}
\biggl[\frac{\sin(\omega_i-\omega_j)}{\omega_i-\omega_j}
+ \frac{\sin(\omega_i+\omega_j)}{\omega_i+\omega_j}
\biggr]_{i,j=0,\ldots,q}.
\]
\end{exmp}
\begin{exmp} Consider now the case $P=Q=1$ in model (\ref{ext-ARFIMA}) where
$\sigma(u)=1$ and $d(u)$, $\Phi(u,B)$, $\Theta(u,B)$ are specified by
\begin{eqnarray*}
d(u) & = & \alpha_1 u, \\
\Phi(u,B) & = & 1+\phi(u) B,\qquad \phi(u)=\alpha_2u,\\
\Theta(u,B) & = & 1+\theta(u) B, \qquad \theta(u)=\alpha_3u,
\end{eqnarray*}
for $u\in[0,1]$. In this case, the parameter vector is $\theta
=(\alpha_1,\alpha_2,\alpha_3)'$, with $0<\alpha_1< \frac12$,
$|\alpha_j|<1$, $j=1,2$
and the matrix $\Gamma$ from (\ref{Gamma}) can be written as
\[
\Gamma= \pmatrix{
\gamma_{11} & \gamma_{12} & \gamma_{13} \cr
\gamma_{21} & \gamma_{22} & \gamma_{23} \cr
\gamma_{31} & \gamma_{32} & \gamma_{33}},
\]
where
\begin{eqnarray*}
\gamma_{11} &=& \frac{1}{2\alpha_1^3}\log\frac{1+\alpha_1}{1-\alpha
_1} -
\frac{1}{\alpha_1^2},\qquad
\gamma_{12} = \frac{g(\alpha_1\alpha_2)}{(\alpha_1\alpha_2)^{3/2}}
- \frac{1}{\alpha_1\alpha_2},\\
\gamma_{13} &=& \frac{1}{2\alpha_1} \biggl\{ \biggl[\frac{1}{2} -
\frac{1}{\alpha_1} \biggr] - \biggl[1- \frac{1}{\alpha_1^2}
\biggr]\log(1 + \alpha_1) \biggr\},\\
\gamma_{22} &=& \frac{1}{2\alpha_2^3}\log\frac{1+\alpha
_2}{1-\alpha_2} -
\frac{1}{\alpha_2^2},\\
\gamma_{23} &=& \frac{1}{2\alpha_2} \biggl\{ \biggl[1- \frac{1}{\alpha_2^2}
\biggr]\log(1 + \alpha_2) - \biggl[\frac{1}{2} -
\frac{1}{\alpha_2} \biggr] \biggr\}, \qquad\gamma_{33} = \frac{\pi^2}{18},
\end{eqnarray*}
with $g(x)=\arctanh(\sqrt{x})$ for $x\in(0,1)$ and $g(x)=\arctan
(\sqrt{-x})$ for
$x\in(-1,0)$.
\end{exmp}

\section{Proofs} \label{proofs}
This section is devoted to the proof of Theorems \ref{CON}--\ref{average-d}. Before
presenting the proofs of these results, we introduce and prove three
useful propositions which are of independent interest. These
propositions involve the large sample properties of the functional
operator defined next. Consider the function $\phi\dvtx
[0,1]\times[-\pi,\pi]\to\mathbb{R} $ and define the functional operator
%
\begin{equation} \label{J}
J(\phi)= \int_0^1 \int_{-\pi}^{\pi} \phi(u,\lambda) f(u,\lambda
) \,d\lambda \,du,
\end{equation}
where $f(u,\lambda)$ is the time-varying spectral density of the limit
process (1). Define
the sample version of $J(\cdot)$ as
%
\begin{equation} \label{JT}
J_T(\phi)= \frac1M \sum_{j=1}^M \int_{-\pi}^{\pi} \phi
(u_j,\lambda) I_N(u_j,\lambda)
\,d\lambda,
\end{equation}
where $M$ and $u_j$, $j=1,\ldots,M$ are given in Section \ref{methodology}.
Furthermore, define the
matrix
%
\begin{equation} \label{Q}
Q(u) = \biggl( \int_{-\pi}^{\pi}\phi(u, \lambda) e^{i\lambda(s-t)
\,d\lambda}
\biggr)_{s,t=1, \ldots, N},
\end{equation}
and the block-diagonal matrix $Q(\phi)=\operatorname{diag}[Q(u_1), \ldots,
Q(u_M)]$. For notational simplicity, sometimes in what follows we have
dropped $\theta$ from $d_{\theta}(u)$ so that it becomes $d(u)$.
\begin{rem}
Since the function $A(u,\lambda)$ and the spectral density
$f(u,\lambda)$ of a locally
stationary long-memory process are unbounded at zero frequency, the
techniques used next
to prove the large sample properties of $J(\phi)$ and the
quasi-likelihood estimators
are different from those used in the short-memory context. For
instance, the function
$A(u,\lambda)$ does not satisfy the key assumption A.1 of \citet
{Dahl97} or the
coefficients $\psi_j( \frac{t}{T})$ of (\ref{wold-general}) fail to meet
conditions (2) and (3) of \citet{Dahl09}. Due to the unboundeness of
$f(u,\lambda)$ at the
origin, our proofs exploit the properties of the Fourier transforms
\begin{eqnarray*}
\widehat{f}(u,\cdot)&=&\int_{-\pi}^{\pi} f(u,\lambda) e^{i \lambda
\cdot} \,d\lambda,\\
\widehat{f}(u,v,\cdot):\!&=&\int_{-\pi}^{\pi} A(u,\lambda
)A(v,-\lambda) e^{i \lambda\cdot} \,d\lambda.
\end{eqnarray*}
\end{rem}

\subsection{Propositions}

\begin{proposition} \label{E}
Let $f(u,\lambda)$ be a time-varying spectral density satisfying
assumption \textup{\hyperlink{A1}{A1}}
and assume that the function $\phi(u,\lambda)$ appearing in (\ref
{J}) is symmetric in
$\lambda$ and twice differentiable with respect to $u$. Let $\widehat
{f}(u, k)$ and
$\widehat{\phi}(u,k)$ be their Fourier coefficients, respectively. If
there is a
positive constant $K$ such that
\[
|\widehat{f}(u,k)\widehat{\phi}(u,k)| \leq K \biggl(\frac{\log k }{k^2} \biggr),
\]
for all $u\in[0,1]$ and $k>1$, then, under assumptions \textup{\hyperlink{A3}{A2}} and
\textup{\hyperlink{A4}{A3}} we have
that
\[
\esp[J_T(\phi)] = J(\phi) + \mathcal{O} \biggl(\frac{\log^2 N}{N}
\biggr) + \mathcal{O} \biggl(\frac{1}{M} \biggr).
\]
\end{proposition}
\begin{pf}
From  definition (\ref{JT}), we can write
\begin{eqnarray*}
\esp[J_T(\phi)] &=& \frac{1}{M} \sum_{j=1}^M \int_{-\pi}^{\pi}
\phi(u_j,
\lambda) \esp[I_N(u_j, \lambda)] \,d \lambda\\
&=& \frac{1}{2\pi M H_{2,N}(0)} \sum_{j=1}^M \int_{-\pi}^{\pi}
\phi(u_j, \lambda) \esp|D_N (u_j, \lambda)|^2 \,d\lambda\\
&=& \frac{1}{2\pi M H_{2,N}(0)} \sum_{j=1}^M \int_{-\pi}^{\pi} \phi
(u_j, \lambda)
\sum_{t, s=0}^{N-1} h \biggl( \frac{t}{N} \biggr) h
\biggl( \frac{s}{N} \biggr)\\
&&\qquad\quad\hspace*{126.6pt}{}\times c(u_j,t,s)e^{i\lambda(s-t)} \,d\lambda,
\end{eqnarray*}
where
\[
c(u,t,s) = \esp\bigl(Y_{[uT]-N/2+t+1,T} Y_{[uT]-N/2+s+1,T} \bigr).
\]
Thus,
\begin{eqnarray*}
\esp[J_T(\phi)]
&=& \frac{1}{2\pi M H_{2,N}(0)}\\
&&\hspace*{0pt}{}\times  \sum_{j=1}^{M}
\sum_{t, s=0}^{N-1} h \biggl( \frac{t}{N} \biggr)
h \biggl( \frac{s}{N} \biggr) c(u_j,t,s) \int_{-\pi}^{\pi} \phi(u_j,
\lambda)e^{i\lambda(s-t)} \,d\lambda\\
&=& \frac{1}{2\pi M H_{2,N}(0)} \sum_{j=1}^{M} \sum_{t, s=0}^{N-1} h
\biggl( \frac{t}{N} \biggr)h
\biggl( \frac{s}{N} \biggr) c(u_j,t,s)\widehat{\phi}(u_j,s-t) \\
&=& \frac{1}{2\pi M H_{2,N}(0)}\sum_{j=1}^{M} \sum_{t=0}^{N-1} \sum
_{k=0}^{N-t}
h \biggl( \frac{t}{N} \biggr)h \biggl( \frac{t}{N}+\frac{k}{N} \biggr)
 c(u_j,t,t+k) \\
&&\hspace*{116.1pt}{}\times \widehat{\phi}(u_j,k)(2-\delta_k),
\end{eqnarray*}
where $\delta_k=1$ for $k=0$ and $\delta_k=0$ for $k\neq0$. By
assumption \hyperlink{A3}{A2}
and Taylor's theorem,
\[
h \biggl( \frac{t}{N} + \frac{k}{N} \biggr) = h \biggl( \frac{t}{N}
\biggr)+ h'(\xi_{t,k,N}) \frac{k}{N},
\]
for some $\xi_{t,k,N} \in( \frac{t}{N}, \frac{t+k}{N} )$,
for $k\geq0$. Thus,
%
\begin{eqnarray} \label{0}\qquad
\esp[J_T(\phi)] &=& \frac{1}{2\pi M
H_{2,N}(0)} \sum_{j=1}^M\sum_{t=0}^{N-1}\sum_{k=0}^{N-t}
h^2 \biggl( \frac{t}{N} \biggr) c(u_j, t, t+k) \nonumber\\
&&\hspace*{128.8pt}\hspace*{-13pt}{}\times\widehat{\phi}(u_j,k) (2-\delta
_k) \nonumber\\[-8pt]\\[-8pt]
&&{} + \frac{1}{2\pi M H_{2,N}(0)} \sum_{j=1}^M\sum_{t=0}^{N-1}\sum_{k=0}^{N-t}
h \biggl( \frac{t}{N} \biggr) h' (\xi_{t,k,N}) c(u_j, t, t+k)\nonumber\\
&&\hspace*{128.8pt}{}\times
\widehat{\phi}(u_j,k) (2-\delta_k).\nonumber
\end{eqnarray}
Under assumption \hyperlink{A1}{A1}, we can expand $c(u, t, t+k)$ by Taylor's
theorem as
\begin{eqnarray*}
&&
c(u,t,t+k) \\
&&\qquad= \widehat{f}(u,k) + \widehat{f}(u,k)\varphi_1(u,k)
\biggl( \frac{t+1-N/2}{T} \biggr) \\
&&\qquad\quad{} + \widehat{f}(u,k) \varphi_2(u,k) \biggl( \frac{k}{T} \biggr) +
R(u,t,k,N,T),
\end{eqnarray*}
where
\begin{eqnarray*}
\varphi_1(u,k) &=& \frac{C_1(\theta,u,u)}{C(\theta,u,u)}+2 d'(u)\log
k,\\
 C_1(\theta,u,u)&=&\frac{\partial C(\theta,u,u+v)}{\partial u}
\bigg|_{v=0},\\
\varphi_2(u,k) &=& \frac{C_2(\theta,u,u)}{C(\theta,u,u)}+d'(u)\log k,
\\
C_2(\theta,u,u)&=&\frac{\partial C(\theta,u,u+v)}{\partial
v} \bigg|_{v=0},
\end{eqnarray*}
$d'(u)=\frac{\partial d_{\theta}(u)}{\partial u}$, $C(\theta,u,v)$
is defined in
assumption \hyperlink{A1}{A1} and the remainder term is given by
\[
R(u,t,k,N,T) = \mathcal{O} \biggl\{ \widehat{f}(u,k) \biggl[
\biggl( \frac{k}{T} \biggr)^2 + \biggl( \frac{t}{T} \biggr)^2 \biggr]
\log^2 k \biggr\}.
\]
Thus, since by assumption \hyperlink{A1}{A1} $|d'(u)|\leq K$ for all $u\in
[0,1]$, we have $|\varphi_j(u,k)|\leq K\log k $ for $j=1,2$ and $k>1$.
Now we can write
%
\begin{eqnarray}\label{(*)}
&&\sum_{k=0}^{N-t} c(u_j, t, t+k) \widehat{\phi}(u_j, k) (2-\delta
_k)\nonumber\\
&&\qquad= \sum_{k=0}^{N-t} \widehat{f}(u_j, k)\widehat{\phi}
(u_j,
k)(2-\delta_k) \nonumber\\
&&\qquad\quad{} + \sum_{k=0}^{N-t} \widehat{f}(u_j, k)\widehat{\phi} (u_j, k)
\varphi_1 (u_j, k)
(2- \delta_k) \biggl( \frac{t+1-N/2}{T} \biggr) \\
&&\qquad\quad{} + \sum_{k=0}^{N-t} \widehat{f}(u_j, k)\widehat{\phi} (u_j, k)
\varphi_2 (_j, k) (2-
\delta_k) { \frac{k}{T}}\nonumber\\
&&\qquad\quad{} + \sum_{k=0}^{N-t} R(u_j, t, k, N, T)
\widehat{\phi} (u_j, k) (2- \delta_k). \nonumber
\end{eqnarray}
Since by assumption $|\widehat{f}(u,k) \widehat{\phi}(u,k)| \leq K
\log k/k^2$, for
$k>1$, uniformly in $u \in[0,1]$, we conclude that there is a finite limit
$A(u)< \infty$ such
that
\[
A(u) = \lim_{N\to\infty} \sum_{k=0}^{N} \widehat{f}(u,k) \widehat
{\phi}(u,k)
(2-\delta_k).
\]
Consequently,
\begin{eqnarray*}
&& \sum_{t=0}^{N-1} h^2 \biggl( \frac{t}{N} \biggr) \sum_{k=0}^{N-t}
\widehat{f}(u_j,k) \widehat{\phi}(u_j, k) (2-\delta_k) \\
&&\qquad= \sum_{t=0}^{N-1} h^2 \biggl( \frac{t}{N} \biggr) \sum_{k=0}^{N-1}
\widehat{f}(u_j, k) \widehat{\phi}(u_j, k) (2-\delta_k) \\
&&\qquad\quad{} - \sum_{t=0}^{N-1} h^2 \biggl( \frac{t}{N} \biggr) \sum_{k=N-t+1}^{N-1}
\widehat{f}(u_j, k) \widehat{\phi}(u_j, k) (2-\delta_k) \\
&&\qquad= \sum_{t=0}^{N-1} h^2 \biggl( \frac{t}{N} \biggr) \Biggl[A(u_j) -
\sum_{k=N}^{\infty} \widehat{f}(u_j, k) \widehat{\phi}(u_j, k)
(2-\delta_k) \Biggr] \\
&&\qquad\quad{}+
\mathcal{O} (\log^2N),
\end{eqnarray*}
by Lemma \ref{Lemma 101}. Hence,
\begin{eqnarray*}
&& \sum_{t=0}^{N-1} h^2 \biggl( \frac{t}{N} \biggr) \sum_{k=0}^{N-t}
\widehat{f}(u_j, k) \widehat{\phi}(u_j, k) (2-\delta_k)\\
&&\qquad = A(u_j)
\sum_{t=0}^{N-1} h^2 \biggl( \frac{t}{N} \biggr)
- \sum_{t=0}^{N-1} h^2 \biggl( \frac{t}{N} \biggr) \sum_{k=N}^{\infty}
\widehat{f}(u_j, k) \widehat{\phi}(u_j, k) (2-\delta_k)\\
&&\qquad\quad{} + \mathcal
{O} (\log^2 N).
\end{eqnarray*}
But,
\[
\Biggl| \sum_{k=N}^{\infty}\widehat{f}(u_j, k) \widehat{\phi}(u_j, k)
(2-\delta_k) \Biggr| < K \sum_{k=N}^{\infty} \frac{\log k}{k^2} =
\mathcal{O} \biggl(\frac{\log N}{N} \biggr),
\]
and consequently,
\[
\Biggl| \sum_{t=0}^{N-1} h^2 \biggl( \frac{t}{N} \biggr) \sum_{k=N}^{\infty}
\widehat{f}(u_j, k) \widehat{\phi}(u_j, k) (2-\delta_k) \Biggr| =
\mathcal{O} (\log N).
\]
Therefore,
\[
\sum_{t=0}^{N-1} h^2 \biggl( \frac{t}{N} \biggr) \sum_{k=0}^{N-t}
\widehat{f}(u_j, k) \widehat{\phi}(u_j, k) (2-\delta_k) = A(u_j)
\sum_{t=0}^{N-1} h^2
\biggl( \frac{t}{N} \biggr) + \mathcal{O} (\log^2 N).
\]
On the other hand, by analyzing the term involving the second summand
of (\ref{(*)}) we
get
\begin{eqnarray*}
&& \sum_{t=0}^{N-1} h^2 \biggl( \frac{t}{N} \biggr)
\sum_{k=0}^{N-t} \varphi_1 (u_j, k) \widehat{f}(u_j,k) \widehat
{\phi}(u_j, k)
(2-\delta_k)
\biggl( \frac{t+1 - N/2}{T} \biggr) \\
&&\qquad = \sum_{t=0}^{N-1} h^2 \biggl( \frac{t}{N} \biggr) \Biggl[ \sum_{k=0}^{N-1}
\varphi_1 (u_j, k) \widehat{f}(u_j,k) \widehat{\phi}(u_j,
k) (2-\delta_k) \Biggr] \biggl(\frac{t+1 - N/2}{T} \biggr) \\
&&\qquad\quad{} + \sum_{t=0}^{N-1} h^2 \biggl( \frac{t}{N} \biggr) \sum_{k=N-t+1}^{N-1}
\varphi_1 (u_j, k) \widehat{\phi}(u_j, k) \widehat{f}(u_j,k)
\biggl( \frac{t+1 - N/2}{T} \biggr)
\\
&&\qquad = \Biggl[ \sum_{t=0}^{N-1} h^2 \biggl( \frac{t}{N} \biggr)
\biggl( \frac{t+1 - N/2}{T} \biggr) \Biggr] \Biggl[ \sum_{k=0}^{N-1} \varphi_1
(u_j, k) \widehat{f}(u_j,k) \widehat{\phi}(u_j, k) (2-\delta_k) \Biggr]
\\
&&\qquad\quad{} + \mathcal{O} \biggl(\frac{N \log^2 N}{T} \biggr),
\end{eqnarray*}
by Lemma \ref{Lemma 102}. Now, since $h(\cdot)$ is symmetric around
$1/2$, we have
\[
\sum_{t=0}^{N-1} h^2 \biggl( \frac{t}{N} \biggr) \biggl( \frac{t+1 -
N/2}{T} \biggr)= \mathcal{O} \biggl(\frac{1}{T} \biggr).
\]
Besides, $| {\sum_{k=0}^{N-1}} \varphi_1(u_j, k) \widehat{f}(u_j, k)
\widehat{\phi}(u_j,
k) (2-\delta_k)| \leq K \sum_{k=1}^{N} \frac{(\log k)^2}{k^2} <
\infty$. Consequently,
\begin{eqnarray*}
&&\sum_{t=0}^{N-1} h^2 \biggl( \frac{t}{N} \biggr) \sum_{k=0}^{N-t} \varphi_1
(u_j, k) \widehat{f}(u_j,k) \widehat{\phi}(u_j, k) (2-\delta_k)
\biggl( \frac{t+1 - N/2}{T} \biggr)\\
&&\qquad =\mathcal{O} \biggl(\frac{N}{T} \log^2
N \biggr).
\end{eqnarray*}
The third term of (\ref{(*)}) can be bounded as follows:
\[
\Biggl|\sum_{k=0}^{N-t} \varphi_2(u_j, k) \widehat{f}(u_j, k) \widehat
{\phi} (u_j, k)
k \Biggr| \leq K \sum_{k=1}^{N} \frac{\log k}{k} \leq K \log^2 N,
\]
and then
\[
\Biggl|\sum_{t=0}^{N-1} h^2 \biggl( \frac{t}{N} \biggr) \sum_{k=0}^{N-t}
\varphi_2(u_j, k) \widehat{f}(u_j, k) \widehat{\phi} (u_j, k) \frac{k}{T}
\Biggr| \leq K \frac{N}{T} \log^2 N.
\]
The last term of (\ref{(*)}) can be bounded as follows:
\[
\Biggl|\sum_{k=0}^{N-t} R(u_j, t, k, N, T) \widehat{\phi}(u_j, k)
(2-\delta_k) \Biggr|
\leq K \log^2 N \biggl(\frac{N}{T} \biggr)^2,
\]
and then
\[
\Biggl|\sum_{t=0}^{N-1} \sum_{k=0}^{N-t} h^2 \biggl( \frac{t}{N} \biggr) R(u_j,
t, k, N, T) \widehat{\phi}(u_j, k) (2-\delta_k) \Biggr| \leq K \frac
{N^3}{T^2} \log^2 N.
\]
Note that by assumption \hyperlink{A4}{A3}, the term above converges to zero as
$N,T\to\infty$.
Therefore, the first term in (\ref{0}) can be written as
\begin{eqnarray*}
&&\frac{1}{2\pi M H_{2,N}(0)} \sum_{j=1}^M\sum_{t, k=0}^{N-1} h^2
\biggl( \frac{t}{N} \biggr) c(u_j,t,t+h) \widehat{\phi} (u_j,h) (2-\delta_k)\\
&&\qquad=
\frac{1}{2 \pi M} \sum_{j=1}^{M} A(u_j) + \mathcal{O} \biggl(\frac{\log
^2 N}{N} \biggr).
\end{eqnarray*}
Now, by Lemma \ref{Lemma A.1} we can write $A(u) = 2\pi\int
_{-\pi}^{\pi} \phi(u,
\omega) f(u, \omega) \,d\omega$, and then
%
\begin{equation} \label{(10)}
\frac{1}{2 \pi M} \sum_{j=1}^{M} A(u_j) = J(\phi) + \mathcal{O}
\biggl(\frac{1}{M} \biggr).
\end{equation}
On the other hand, the second term in (\ref{0}) can be bounded as
follows:
\begin{eqnarray*}
&&|c(u_j, t, t+k) \widehat{\phi}(u_j, k) (2-\delta_k)| \\
&&\qquad\leq K \biggl\{
\widehat{f}(u_j, k) \widehat{\phi} (u_j, k) + \frac{N}{T}|\widehat{f}(u_j, k) \widehat{\phi} (u_j, k) \varphi_1
(u_j,k)| \\
&&\qquad\quad\hspace*{14pt}{} + \frac{N}{T}|\widehat{f}(u_j, k) \widehat{\phi} (u_j, k) \varphi
_2 (u_j, k)| +
\frac{N^2}{T^2} \log^2 N |\widehat{f}(u_j, k) \widehat{\phi} (u_j,
k)| \biggr\}.
\end{eqnarray*}
Since $|\varphi_i(u_j, k) | \leq K \log k$ for $i=1,2, j=1, \ldots,
M$ and $k>1$, we
conclude that
\[
|c(u_j, t, t+k) \widehat{\phi}(u_j, k) (2-\delta_k)| \leq
K \frac{N}{T} \frac{\log k}{k^2}.
\]
Therefore, since $|h'(u)| \leq K$ for $u \in[0,1]$ by assumption
\hyperlink{A3}{A2}, we have
\[
\Biggl|\sum_{k=0}^{N-t} c(u_j, t, t+k) \widehat{\phi}(u_j, k) \frac{k}{N}
h'(\xi_{t,k,N}) \Biggr| \leq\frac{K}{T} \sum_{k=1}^{N} \frac{\log k}{k}
\leq K
\frac{\log^2 N}{T}.
\]
Consequently,
\[
\Biggl|\sum_{t=0}^{N-1} h \biggl( \frac{t}{N} \biggr) \sum_{k=0}^{N-t} c(u_j, t,
t+k) \widehat{\phi}(u_j, k) (2-\delta_k) \frac{k}{N} h' (\xi
_{t,k,N}) \Biggr| \leq K N
\frac{\log^2 N}{T}.
\]
Hence, the second term of (\ref{0}) is bounded by $K (\log^2 N)/T$.
From this and
(\ref{(10)}), the required result is obtained.
\end{pf}
%
%
\begin{proposition} \label{Cov} 
Let $f(u,\lambda)$ be a time-varying spectral density satisfying
assumption \textup{\hyperlink{A1}{A1}}.
Let $\phi_1, \phi_2\dvtx [0,1]\times[-\pi,\pi]\to\mathbb{R} $ be two
functions such that
$\phi_1(u,\lambda)$ and $\phi_2(u,\lambda)$ are symmetric in
$\lambda$, twice
differentiable with respect to $u$ and their Fourier coefficients satisfy
$|\widehat{\phi}_1(u,k)|,|\widehat{\phi}_2(u,k)| \leq K
|k|^{-2d(u)-1}$ for
$u\in[0,1]$ and $|k|>1$. If assumptions \textup{\hyperlink{A3}{A2}} and
\textup{\hyperlink{A4}{A3}}
hold, then
\[
\lim_{T\to\infty} T \cov[J_T(\phi_1), J_T(\phi_2) ] = 4\pi
\int_0^1\int_{-\pi}^{\pi}\phi_1(u,\lambda)\phi_2(u,\lambda)
f(u,\lambda)^2
\,d\lambda \,du.
\]
\end{proposition}
\begin{pf} We can write
\begin{eqnarray*}
&&
T \cov[J_T(\phi_1), J_T(\phi_2)] \\
&&\qquad = \frac{T}{M^2} \sum_{j,k = 1}^{M}
\int_{-\pi}^{\pi}\int_{-\pi}^{\pi} \phi_1(u_j,\lambda) \phi
_2(u_k, \mu)\\
&&\qquad\quad\hspace*{79.3pt}{}\times\cov[I_N(u_j,\lambda), I_N(u_k,\mu)] \,d\lambda \,d\mu.
\end{eqnarray*}
But,
\begin{eqnarray*}
&&\cov[I_N(u_j,\lambda), I_N(u_k,\mu)] \\
&&\qquad= \frac{1}{[2\pi
H_{2,N}(0)]^2}\cov(|D_N(u_j,\lambda)|^2, D_N(u_k,\mu)|^2)\\
&&\qquad = \frac{1}{[2\pi H_{2,N}(0)]^2} \sum_{t,s,p,m=0}^{N-1}
h \biggl( \frac{t}{N} \biggr)h \biggl( \frac{s}{N} \biggr)h \biggl( \frac{p}{N} \biggr)h \biggl( \frac{m}{N}
\biggr)\\
&&\qquad\quad\hspace*{102pt}{}\times
e^{\imag\lambda(s-t) + \imag\mu(m-p)}\\
&&\qquad\quad\hspace*{102pt}{} \times\cov\bigl(Y_{[u_jT]-{N/2} + s + 1,T}
Y_{[u_jT]-{N/2} + t + 1,T},\\
&&\qquad\quad\hspace*{135pt}Y_{[u_kT]-{N/2} + p + 1,T}
Y_{[u_kT]-{N/2} + m + 1,T} \bigr).
\end{eqnarray*}
Now, an application of Theorem 2.3.2 of \citet{Bril81} yields
\begin{eqnarray*}
&&\cov[I_N(u_j,\lambda), I_N(u_k,\mu)]
\\
&&\qquad= \frac{1}{[2\pi H_{2,N}(0)]^2}\\
&&\qquad\quad{}\times \sum_{t,s,p,m=0}^{N-1} h \biggl( \frac
{t}{N} \biggr)
h \biggl( \frac{s}{N} \biggr)h \biggl( \frac{p}{N} \biggr)
h \biggl( \frac{m}{N} \biggr)  e^{\imag\lambda(s-t) + \imag\mu(m-p)}\\
&&\qquad\quad\hspace*{48.7pt}{} \times
\bigl\{ \cov\bigl(Y_{[u_jT]-{N/2} + t + 1,T}, Y_{[u_k
T]-{N/2} + m + 1,T} \bigr) \\
&&\qquad\quad\hspace*{63.7pt}{} \times\cov\bigl(Y_{[u_jT]-{N/2} + s + 1,T},
Y_{[u_k T]-{N/2} + p + 1,T}\bigr)\\
&&\qquad\quad\hspace*{63.7pt}{} + \cov\bigl(Y_{[u_jT]-{N/2} + t + 1,T}, Y_{[u_k T]-{N/2} + p
+ 1,T}\bigr)\\
&&\qquad\quad\hspace*{73.6pt}{} \times\cov\bigl(Y_{[u_jT]-{N/2} + s + 1,T}, Y_{[u_k T]-
{N/2} + p + 1,T} \bigr) \bigr\}\\
&&\qquad = \frac{1}{[2\pi H_{2,N}(0)]^2} \\
&&\qquad\quad{}\times\int_{-\pi}^{\pi}\int_{-\pi
}^{\pi}
H_N \biggl(A^0_{t_j-{N/2}+1+\cdot, T}(x)
h \biggl( \frac{\cdot}{N} \biggr),\lambda-x \biggr)\\
&&\qquad\quad\hspace*{49.3pt}{} \times H_N \biggl(A^0_{t_k-{N/2}+1+\cdot, T}(x)
h \biggl( \frac{\cdot}{N} \biggr),x-\mu\biggr)\\
&&\qquad\quad\hspace*{49.3pt}{}\times H_N \biggl(A^0_{t_j-{N/2}+1+\cdot, T}(y) h \biggl( \frac{\cdot
}{N} \biggr),-y-\lambda\biggr) \\
&&\qquad\quad\hspace*{49.3pt}{} \times H_N \biggl(A^0_{t_k-{N/2}+1+\cdot, T}(y)
h \biggl( \frac{\cdot}{N} \biggr),y+\mu\biggr)\\
&&\qquad\quad\hspace*{49.3pt}{} \times e^{\imag\lambda(s-t) + \imag
\mu(m-p)} \,dx \,dy
\\
&&\qquad\quad{} + \frac{1}{[2\pi H_{2,N}(0)]^2} \int_{-\pi}^{\pi}\int_{-\pi
}^{\pi} H_N \biggl(A^0_{t_j-{N/2}+1+\cdot, T}(x) h \biggl( \frac{\cdot
}{N} \biggr),\lambda-x \biggr)\\
&&\qquad\quad\hspace*{116pt}{} \times H_N \biggl(A^0_{t_k-{N/2}+1+\cdot, T}(x)
h \biggl( \frac{\cdot}{N} \biggr),x+\mu\biggr)\\
&&\qquad\quad\hspace*{116pt}{} \times
H_N \biggl(A^0_{t_j-{N/2}+1+\cdot, T}(y) h \biggl( \frac{\cdot}{N}
\biggr),-y-\lambda\biggr) \\
&&\qquad\quad\hspace*{116pt}{} \times  H_N \biggl(A^0_{t_k-{N/2}+1+\cdot, T}(y)
h \biggl( \frac{\cdot}{N} \biggr),y-\mu\biggr) \\
&&\qquad\quad\hspace*{116pt}{} \times e^{\imag\lambda(s-t) + \imag
\mu(m-p)} \,dx \,dy.
\end{eqnarray*}
Thus,
%
\begin{equation} \label{Tcov}
T\cov(J_T(\phi_1), J_T(\phi_2)) = \frac{T}{ [2\pi M H_{2,N}(0) ]^2}
\bigl[B_N^{(1)} + B_N^{(2)} \bigr],
\end{equation}
where 
\begin{eqnarray*}
B_N^{(1)} &=& \sum_{j,k = 1}^{M} \int_{\Pi} \phi_1(u_j,
\lambda)\phi_2(u_k,\mu) H_N \biggl(A_{t_j-{N/2}+1+\cdot,T}^{0} (x)
h \biggl( \frac{\cdot}{N} \biggr), \lambda- x \biggr) \\
&&\hspace*{34.64pt}{} \times H_N \biggl(\overline{A_{t_k-{N/2}+1+\cdot,T}^{0} (x)} h \biggl(
\frac{\cdot}{N} \biggr), x-\mu\biggr) H_N \\
&&\hspace*{34.64pt}{} \times \biggl(A_{t_j-{N/2}+1+\cdot,T}^{0}
(y) h \biggl( \frac{\cdot}{N} \biggr), -y-\lambda\biggr) \\
&&\hspace*{34.64pt}{} \times H_N \biggl(\overline{A_{t_k-{N/2}+1+\cdot,T}^{0} (y)} h \biggl(
\frac{\cdot}{N} \biggr), y+\mu\biggr) \\
&&\hspace*{34.64pt}{} \times e^{\imag(x+y)(t_j-t_k)} \,dx \,dy \,d\mu
\,d\lambda,
\end{eqnarray*}
with $\Pi= [-\pi,\pi]^4$, and
\begin{eqnarray*}
B_N^{(2)} &=& \sum_{j,k = 1}^{M} \int_{\Pi} \phi_1(u_j, \lambda
)\phi_2(u_k,\mu)
H_N \biggl(A_{t_j-{N/2}+1+\cdot,T}^{0} (x)
h \biggl( \frac{\cdot}{N} \biggr), \lambda- x \biggr) \\
&&\hspace*{34.9pt}{} \times H_N \biggl(\overline{A_{t_k-{N/2}+1+\cdot,T}^{0} (x)} h \biggl(
\frac{\cdot}{N} \biggr), x+\mu\biggr)\\
&&\hspace*{34.9pt}{} \times  H_N \biggl(A_{t_j-{N/2}+1+\cdot,T}^{0}
(y) h \biggl( \frac{\cdot}{N} \biggr), -y-\lambda\biggr) \\
&&\hspace*{34.9pt}{} \times  H_N \biggl(\overline{A_{t_k-{N/2}+1+\cdot,T}^{0} (y)} h \biggl(
\frac{\cdot}{N} \biggr), y-\mu\biggr)\\
&&\hspace*{34.9pt}{} \times  e^{\imag(x+y)(t_j-t_k)} \,dx \,dy \,d\mu
\,d\lambda.
\end{eqnarray*}
The term $B_N^{(1)}$ can be written as follows: 
%
\begin{eqnarray} \label{BN}\hspace*{18pt}
B_N^{(1)} &=& \sum_{j,k = 1}^{M}
\int_{\Pi} \phi_1(u_j,
\lambda)\phi_2(u_k,\mu)\nonumber\\
&&\hspace*{34.6pt}{} \times A(u_j, x) A(u_k, -x) A(u_j, y) A(u_k, -y) \nonumber\\
&&\hspace*{34.6pt}{} \times H_N(\lambda-x) H_N(x-\mu) H_N(\mu+y) H_N(-y-\lambda)
\nonumber\\
&&\hspace*{34.6pt}{}\times e^{\imag(x+y)(t_j-t_k)} \,dx \,dy \,d\lambda \,d\mu
+ R_N \nonumber\\[-8pt]\\[-8pt]
& = &\sum_{j,k = 1}^{M} \int_{\Pi}
\phi_1(u_j, x) A(u_j, x) A(u_k, -x)\phi_2(u_k,y) A(u_j, y)\nonumber\\
&&\hspace*{54.1pt}{} \times A(u_k, -y)
 H_N(\lambda-x) H_N(x-\mu) H_N(\mu+y)\nonumber\\
&&\hspace*{54.1pt}{}\times H_N(-y-\lambda)
e^{\imag(x+y)(t_j-t_k)} \,dx \,dy \,d\lambda \,d\mu\nonumber\\
&&{} + \Delta_N^{(1)} + \Delta_N^{(2)} + R_N,\nonumber
\end{eqnarray}
with
\begin{eqnarray*}
\Delta_N^{(1)} &=& \sum_{j,k = 1}^{M} \int_{\Pi}
[\phi_1(u_j, \lambda)-\phi_1(u_j,x) ]\phi_2(u_k,\mu) A(u_j, x)
A(u_k, -x) A(u_j, y)\\
&&\hspace*{33.6pt}{} \times A(u_k, -y) H_N(\lambda-x) H_N(x-\mu) H_N(\mu+y)\\
&&\hspace*{33.6pt}{} \times
H_N(-y-\lambda) e^{\imag(x+y)(t_j-t_k)} \,dx \,dy \,d\lambda \,d\mu,
\\
\Delta_N^{(2)} &=& \sum_{j,k = 1}^{M} \int_{\Pi} \phi_1(u_j,x)
[\phi_2(u_k,
\mu)-\phi_2(u_k,y) ]
A(u_j, x) A(u_k, -x) A(u_j, y)\\
&&\hspace*{34.6pt}{} \times A(u_k, -y) H_N(\lambda-x) H_N(x-\mu) H_N(\mu+y)
\\
&&\hspace*{34.6pt}{}\times H_N(-y-\lambda) e^{\imag
(x+y)(t_j-t_k)} \,dx \,dy \,d\lambda \,d\mu,
\end{eqnarray*}
and by Lemma \ref{A.5} the remainder term $R_N$ can be bounded as
follows:
%
\begin{eqnarray} \label{RN}\qquad
|R_N| &\leq& \frac{N}{T} \Biggl|\sum_{j,k = 1}^{M}
\int_{\Pi} \phi_1(u_j,\lambda)
\phi_2(u_k,\mu) A(u_j,x)
A(u_k,-x) \nonumber\\
&&\hspace*{49.7pt}{} \times A(u_j,y) y^{-d(u_k)}
L_N(y+\mu)\nonumber\\[-8pt]\\[-8pt]
&&\hspace*{49.7pt}{} \times  H_N(\lambda- x) H_N(x-\mu)\nonumber\\
&&\hspace*{49.7pt}{} \times H_N(-y-\lambda) e^{\imag(x+y)(t_j-t_k)} \,dx \,dy\,
d\lambda \,d\mu\Biggr|. \nonumber
\end{eqnarray}
By integrating with respect to $\mu$ the term $\Delta_N^{(1)}$ can be
written as
\begin{eqnarray*}
\Delta_N^{(1)} &=& \sum_{j,k=1}^{M} \sum_{t,s=0}^{N-1} h \biggl( \frac
{t}{N} \biggr) h \biggl( \frac{s}{N} \biggr) \widehat{\phi}_2(u_k,t-s)\\
&&\hspace*{41.8pt}{} \times\int_{-\pi}^{\pi}\int_{-\pi}^{\pi}\int_{-\pi}^{\pi}
[\phi_1(u_j,\lambda) - \phi_1(u_j,x) ] A(u_j,x) A(u_k,-x)\\
&&\hspace*{110.3pt}{} \times A(u_j,y) A(u_k,-y) H_N(\lambda-x) H_N(-y-\lambda)\\
&&\hspace*{110.3pt}{} \times  e^{\imag
(x+y)(t_j-t_k)-\imag xt-\imag y s} \,dx \,dy \,d\lambda,
\end{eqnarray*}
and by integrating with respect to $y$ we get
\begin{eqnarray*}
\Delta_N^{(1)} & = & \sum_{j,k=1}^{M} \sum_{t,s,p=0}^{N-1} h \biggl( \frac
{t}{N} \biggr)h \biggl( \frac{s}{N} \biggr)h \biggl( \frac{p}{N} \biggr)\\
&&\hspace*{50.2pt}{} \times \widehat{\phi
}_2(u_k,t-s)\widehat{f}(u_j,u_k,t_j-t_k-s+p)\\
&&\hspace*{50.2pt}{} \times \int_{-\pi}^{\pi}\int
_{-\pi}^{\pi} [\phi_1(u_j,\lambda) - \phi_1(u_j,x) ]\\
&&\hspace*{98.2pt}{} \times A(u_j,x) A(u_k,-x) H_N(\lambda-x)\\
&&\hspace*{98.2pt}{} \times e^{\imag x(t_j-t_k-t)+\imag
\lambda p} \,dx \,d\lambda
\\
& = &\sum_{j,k=1}^{M} \sum_{t,s,p=0}^{N-1} h \biggl( \frac{t}{N} \biggr)h \biggl( \frac
{s}{N} \biggr)h \biggl( \frac{p}{N} \biggr) \widehat{\phi}_2(u_k,t-s)\\
&&\hspace*{49.8pt}{} \times\widehat{f}(u_j,u_k,t_j-t_k-s+p)\\
&&\hspace*{49.8pt}{} \times \varepsilon_N(u_j,u_k,p,t_j-t_k-t),
\end{eqnarray*}
where $\widehat{f}(u,v,k)$ and $\varepsilon_N(r)$ are given by
$\widehat{f}(u,v,k)=\int_{-\pi}^{\pi} A(u,\lambda)A(v,-\lambda)
e^{\imag\lambda k}\,
d\lambda$, and
\begin{eqnarray*}
\hspace*{-3pt}&&
\varepsilon_N(u_j,u_k,p,r) \\
\hspace*{-3pt}&&\qquad= \sum_{m=0}^{N-1} h \biggl( \frac{m}{N} \biggr)
\widehat{\phi}_1(u_j,p-m)
\int_{-\pi}^{\pi} A(u_j,x) A(u_k, -x) e^{\imag x(r+m)} \,dx\\
\hspace*{-3pt}&&\qquad\quad{} - 2\pi h \biggl(\frac{p}{m} \biggr) \int_{-\pi}^{\pi}
\phi_1(u_j,x) A(u_j,x) A(u_k,x) e^{\imag x(r+p)} \,dx\\
\hspace*{-3pt}&&\qquad = \int_{-\pi}^{\pi}\int_{-\pi}^{\pi} \phi_1(u_j,\lambda)
A(u_j, x) A(u_k, -x) e^{\imag(p\lambda+r x)}
\sum_{m=0}^{N-1} h \biggl( \frac{m}{N} \biggr) e^{\imag m(x-\lambda)} \,d\lambda
\,dx\\
\hspace*{-3pt}&&\qquad\quad{} - 2\pi\int_{-\pi}^{\pi} h \biggl( \frac{m}{N} \biggr)
\phi_1(u_j,x) A(u_j,x) A(u_k,-x) e^{\imag(r+p)x} \,dx.
\end{eqnarray*}
But $h (\frac{m}{M} ) = h (\frac{p}{M} ) +
h^{\prime} (\xi_{p,m} )\frac{m-p}{N}$ for some
$\xi_{p,m}\in[0,1]$. Thus,
\begin{eqnarray*}
&&\varepsilon_N(u_j,u_k,p,r)\\
&&\qquad = h \biggl( \frac{p}{N} \biggr) \Biggl\{\int_{-\pi}^{\pi
}\int_{-\pi}^{\pi} \phi_1(u_j,\lambda) A(u_j,x) A(u_k,-x) e^{\imag
(p\lambda+r x)}\\
&&\qquad\quad\hspace*{74.2pt}{} \times\sum_{m=0}^{N-1} e^{\imag m (x-\lambda)} \,d\lambda \,dx\\
&&\qquad\quad\hspace*{35.7pt}{} - 2\pi
\int_{-\pi}^{\pi} \phi_1(u_j,x) A(u_j,x) A(u_k,-x) e^{\imag(r+p)x}\,
dx \Biggr\}\\
&&\qquad\quad{} + \sum_{m=0}^{N-1} h^{\prime} (\xi_{p,m} )\frac{m-p}{N}
\int_{-\pi}^{\pi}\int_{-\pi}^{\pi} \phi_1(u_j,\lambda) A(u_j,
x) A(u_k,-x)\\
&&\qquad\quad\hspace*{138.2pt}{}
\times e^{\imag m(x-\lambda)+\imag(p\lambda+rx)} \,d\lambda \,dx\\
&&\qquad = h \biggl( \frac{p}{N} \biggr)
\varepsilon_N^{(1)}(u_j,u_k,p,r)+\varepsilon_N^{(2)}(u_j,u_k,p,r),
\end{eqnarray*}
where the term $\varepsilon_N^{(1)}(u_j,u_k,p,r)$ is given by
\[
\varepsilon_N^{(1)}(u_j,u_k,p,r) = \int_{-\pi}^{\pi} g(\omega)
e^{\imag r\omega}
\sum_{m=0}^{N-1} e^{\imag m \omega} \,d\omega- 2\pi g(0),
\]
with $ g(\omega) = \int_{-\pi}^{\pi} \phi_1(u_j,\lambda)
A(u_j,\omega+\lambda)
A(u_k,-\omega-\lambda) e^{\imag(p+r)\lambda} \,d\lambda$. Observe
that by Lemma
\ref{Lemma A.1}, for every $u_j,u_k,p,r$, $\varepsilon
_N^{(1)}(u_j,u_k,p,r)\to0$ as
$N\to\infty$, consequently we can write
\begin{eqnarray*}
\varepsilon_N^{(1)}(u_j,u_k,p,r) &=& \int_{-\pi}^{\pi} g(\omega)
\sum_{m=N}^{\infty} e^{\imag(m+r) \omega} \,d\omega\\
&=&
\sum_{m=N}^{\infty} \widehat{\phi}_1(u_j, p-m) \widehat{f}(u_j,u_k,r+m).
\end{eqnarray*}
On the other hand, by assumption \hyperlink{A1}{A1}, $|\widehat{f}(u_j,u_k,
r+m)|\leq K
|r+m|^{d(u_j)+d(u_k)-1}$. Thus, the term $\varepsilon
_N^{(2)}(u_j,u_k,p,r)$ is bounded
by
\begin{eqnarray*}
&& \bigl|\varepsilon_N^{(2)}(u_j,u_k,p,r) \bigr| \\
&&\qquad\leq\frac{K}{N} \sum
_{m=0}^{N-1} |m-p|
\widehat{\phi}_1(u_j, m-p) \widehat{f}(u_j,u_k, r+m)\\
&&\qquad \leq\frac{K}{N} \sum_{m=0}^{N-1} |m-p|^{-2d(u_j)}
|r+m|^{d(u_j)+d(u_k)-1}\\
&&\qquad \leq K \Biggl\{ \sum_{m=0}^{N-1} \biggl| \frac{m}{N}- \frac{p}{N} \biggr|
\biggl| \frac{r}{N}- \frac{m}{N} \biggr|^{d(u_j)+d(u_k)-1} \frac{1}{N} \Biggr\}
N^{d(u_k)-d(u_j)-1}\\
&&\qquad \leq K\int_0^1 \biggl|x- \frac{p}{m} \biggr|^{-2d(u_j)}
\biggl| \frac{r}{N}+x \biggr|^{d(u_j)+d(u_k)-1} \,dx N^{d(u_k)-d(u_j)-1}\\
&&\qquad \leq K N^{d(u_k)-d(u_j)-1} \biggl( \frac{r}{N} \biggr)^{d(u_j)+d(u_k)-1}
\int_0^1 \biggl|x- \frac{p}{m} \biggr|^{-2d(u_j)} \,dx\\
&&\qquad \leq K N^{-2d(u_j)} r^{d(u_j)+d(u_k)-1},
\end{eqnarray*}
where for notational simplicity we have dropped $\theta$ from
$d_{\theta}(\cdot)$. Thus,
\begin{eqnarray*}
\varepsilon_N(u_j,u_k,p,r)& = & h \biggl( \frac{p}{N} \biggr) \sum_{m=N}^{\infty}
\widehat{\phi}_1(u_j,p-m) \widehat{f}(u_j,u_k, r+m) \\
&&{} + \mathcal{O} \bigl(N^{-2d(u_j)} r^{d(u_j)+d(u_k)-1} \bigr).
\end{eqnarray*}
Hence, $\Delta_N^{(1)}$ can be written as
\begin{eqnarray*}
\Delta_N^{(1)} &=& \sum_{j,k=1}^{M} \sum_{t,s,p = 0}^{N-1}
h \biggl( \frac{t}{N} \biggr) h \biggl( \frac{s}{N} \biggr)
h \biggl( \frac{p}{N} \biggr) \widehat{\phi}_2(u_k,t-s) \widehat
{f}(u_j,u_k,t_j-t_k)\\
& &\hspace*{50.1pt}{}\times\Biggl\{h \biggl( \frac{p}{N} \biggr)\sum_{m=N+1}^{\infty}
\widehat{\phi}_1(u_j,p-m) \widehat{f}(u_j,u_k,t_j-t_k+m)
\\
&&\hspace*{117pt}{} + O \bigl(N^{-2d(u_j)} |t_j-t_k-t|^{d(u_j)+d(u_k)-1} \bigr) \Biggr\}\\
&:=&\Delta
_N^{(1.1)} + \Delta_N^{(1.2)},
\end{eqnarray*}
say. Therefore,\vspace*{1pt} $ |\Delta_N^{(1)} | \leq|\Delta_N^{(1.1)} | +
|\Delta_N^{(1.2)} |$. Observe that since $\phi_2(u,\lambda)\sim C
|\lambda|^{2d(u)}$ as $\lambda\to0$ and $d(u)>0$ for all $u\in
[0,1]$, we conclude that
$\phi_2(u,0)={\sum_{k=-\infty}^{\infty}}\widehat{\phi}_2(u$,\break $k)=0$. Thus,
\begin{eqnarray*}
&&
\sum_{t,s = 0}^{N-1}h \biggl( \frac{t}{N} \biggr) h \biggl( \frac{s}{N} \biggr)
\widehat{\phi}_2(u,t-s)\\
&&\qquad=\sum_{t=0}^{N-1} \sum_{k=1-N}^{N-1} h \biggl(
\frac{t}{N} \biggr)
h \biggl( \frac{t}{N}+ \frac{k}{N} \biggr) \widehat{\phi}_2(u,k),
\end{eqnarray*}
where for simplicity we assume that $h(x)=0$ for $x$ outside $[0,1]$.
Now, by an
application of Taylor's theorem we can write
$h ( \frac{t}{N}+ \frac{k}{N} )=h ( \frac{t}{N} )
+h'(\xi(t,k)) \frac{k}{N}$ for some $\xi(t,k) \in
( \frac{t}{N}- \frac{|k|}{N}, \frac{t}{N}+ \frac{|k|}{N} )
$. Hence,
\begin{eqnarray*}
&&\sum_{t=0}^{N-1} \sum_{k=1-N}^{N-1} h \biggl( \frac{t}{N} \biggr)
h \biggl( \frac{t}{N}+ \frac{k}{N} \biggr) \widehat{\phi}_2(u,k)
\\
&&\qquad=
\sum_{t=0}^{N-1} h \biggl( \frac{t}{N} \biggr)^2 \sum_{k=1-N}^{N-1}
\widehat{\phi}_2(u,k)\\
&&\qquad\quad{} +\sum_{t=0}^{N-1} \sum_{k=1-N}^{N-1} h \biggl( \frac{t}{N} \biggr)
h' (\xi(t,k) ) \widehat{\phi}_2(u,k) \frac{k}{N}.
\end{eqnarray*}
Note that $\sum_{k=1-N}^{N-1} \widehat{\phi}_2(u,k)=2\sum
_{k=N}^{\infty}
\widehat{\phi}_2(u,k)$. Therefore, $ |{\sum_{k=1-N}^{N-1}}
\widehat{\phi}_2(u,k) | \leq K \sum_{k=N}^{\infty} k^{-2d(u)-1}\leq
K N^{-2d(u)}$.
Consequently, $ |{\sum_{t=0}^{N-1} h} ( \frac{t}{N} )^2
\sum_{k=1-N}^{N-1} \widehat{\phi}_2(u$,\break $k) |\leq K N^{1-2d(u)}$. On
the other hand,
$ |\sum_{t=0}^{N-1} \sum_{k=1-N}^{N-1} h ( \frac{t}{N} )
h' (\xi(t,k) ) \widehat{\phi}_2(u,k) \times\break\frac{k}{N} |\leq
K \sum_{k=1}^Nk^{-2d(u)}\leq KN^{1-2d(u)}$. Hence,
\[
\Biggl|\sum_{t=0}^{N-1} \sum_{k=1-N}^{N-1} h \biggl( \frac{t}{N} \biggr)
h \biggl( \frac{t}{N}+ \frac{k}{N} \biggr) \widehat{\phi}_2(u,k) \Biggr|\leq
KN^{1-2d(u)}.
\]
Thus, we conclude that
\begin{eqnarray*}
\bigl|\Delta_N^{(1.1)} \bigr| &\leq& K \sum_{j,k=1,j\neq k}^{M} \sum
_{p=0}^{N-1} N^{1-2d(u_k)} |S(j-k)+p|^{d(u_j)+d(u_k)-1}\\
&&\hspace*{69.1pt}{} \times\sum_{m=N+1}^{\infty} |p-m|^{-2d(u_j)-1}
|S(j-k)+m|^{d(u_j)+d(u_k)-1}\\
& \leq & K \sum_{j,k=1,j\neq k}^{M} \sum_{p=0}^{N-1} \biggl|\frac
{S}{N}(j-k)+ \frac{p}{N} \biggr|^{d(u_j)+d(u_k)-1} \\
&&\qquad\quad\hspace*{35.1pt}{}\times\frac{1}{N} \sum
_{m=N+1}^{\infty} \biggl| \frac{p}{N}- \frac{m}{N} \biggr|^{-2d(u_j)-1}\\
&&\qquad\quad\hspace*{90.1pt}{}\times\biggl|\frac{S}{N}(j-k)+ \frac{m}{N} \biggr|^{d(u_j)+d(u_k)-1} \frac
{1}{N}\\
& \leq & K \sum_{j,k=1,j\neq k}^{M} \int_0^1\int_1^{\infty} \biggl|\frac
{S}{N}(j-k)+x \biggr|^{d(u_j)+d(u_k)-1} |x-y |^{-2d(u_j)-1}\\
&&\hspace*{85pt}{} \times\biggl|\frac{S}{N}(j-k)+y \biggr|^{d(u_j)+d(u_k)-1} \,dy \,dx.
\end{eqnarray*}
Let $\delta>0$ and define $I_1(\delta)=\{j,k=1,M\dvtx k<j \vee k-j>
\frac{N}{S}(1+\delta)\}$ and $I_2(\delta)=\{j,k=1,M\dvtx 0<k-j \leq
\frac{N}{S}(1+\delta)\}$.
Therefore, the sum above can be written as $\sum_{j,k=1,j\neq
k}^{M}\cdot=\sum_{I_1(\delta)}\cdot+\sum_{I_2(\delta)}\cdot
:=|\Delta_N^{(1.1.1)} |+ |\Delta_N^{(1.1.2)} |$, say. Observe
that over
$I_1(\delta)$ we have that $ |\frac{S}{N}(j-k)+x |^{-\alpha} \leq K
|\frac{S}{N}(j-k) |^{-\alpha}$ for $\alpha>0$. 
Hence,
\begin{eqnarray*}
\bigl|\Delta_N^{(1.1.1)} \bigr| & \leq & K \sum_{I_1(\delta)} \biggl|\frac{S}{N}(j-k)
\biggr|^{d(u_j)+d(u_k)-1}\\
&&\hspace*{26.7pt}{} \times\int_0^1 \int_1^{\infty} |x-y|^{-2d(u_j)-1}\\
&&\hspace*{74.7pt}{}\times \biggl|\frac
{S}{N}(j-k)+y \biggr|^{d(u_j)+d(u_k)-1}
\,dy \,dx\\
& \leq & K \sum_{j,k=1,j\neq k}^{M} \biggl|\frac{S}{N}(j-k)
\biggr|^{d(u_j)+d(u_k)-1}\\
&&\hspace*{49.3pt}{} \times\int_0^1 \int_1^{\infty} |x-y|^{-2d(u_j)-1}\\
&&\qquad\quad\hspace*{64.4pt}{}\times \biggl|\frac
{S}{N}(j-k)+y \biggr|^{d(u_j)+d(u_k)-1}\,
dy \,dx.
\end{eqnarray*}
Since the integrands in the above expression are all positive, an
application of
Tonelli's theorem yields
\begin{eqnarray*}
\bigl|\Delta_N^{(1.1.1)} \bigr| & \leq & K \sum_{j,k=1,j\neq k}^{M} \biggl|\frac
{S}{N}(j-k) \biggr|^{d(u_j)+d(u_k)-1}\\
&&\hspace*{48.6pt}{} \times\int_1^{\infty} \int_0^1 |x-y|^{-2d(u_j)-1}\\
&&\hspace*{95.9pt}{}\times \biggl|\frac
{S}{N}(j-k)+y \biggr|^{d(u_j)+d(u_k)-1} \,dx \,dy\\
& \leq & K \sum_{j,k=1,j\neq k}^{M} \biggl|\frac{S}{N}(j-k)
\biggr|^{d(u_j)+d(u_k)-1}\\
&&\hspace*{49.3pt}{} \times\int_1^{\infty} \bigl[(y-1)^{-2d(u_j)}-y^{-2d(u_j)} \bigr]\\
&&\hspace*{79.1pt}{}\times \biggl|\frac
{S}{N}(j-k)+y \biggr|^{d(u_j)+d(u_k)-1} \,dy.
\end{eqnarray*}
Then, by Lemma \ref{XXXX} we conclude that
\begin{eqnarray*}
\bigl|\Delta_N^{(1.1.1)} \bigr|&
\leq& K \sum_{j,k=1,j\neq k}^{M} \biggl|\frac{S}{N}(j-k)
\biggr|^{2d(u_j)+2d(u_k)-2}\\
& \leq& K \Biggl[ \sum_{\stackrel{j,k=1,j\neq k}{d(u_j)+d(u_k)\leq1/2}}^M
\biggl|\frac{S}{N}(j-k) \biggr|^{2d(u_j)+2d(u_k)-2} \\
&&\hspace*{13.9pt}{} + \sum_{\stackrel{j,k=1,j\neq
k}{d(u_j)+d(u_k)> 1/2}}^M \biggl|\frac{S}{N}(j-k) \biggr|^{2d(u_j)+2d(u_k)-2} \Biggr].
\end{eqnarray*}
For the first summand above, we have the upper bound
\[
\sum_{j,k=1,j\neq k}^{M} |j-k|^{-1} \biggl(\frac{N}{S} \biggr)^{2} \leq K
\biggl(\frac{N}{S} \biggr)^{2} M \log M,
\]
while the second summand can be bounded as follows:
\begin{eqnarray*}
&& \sum_{\stackrel{j,k=1,j\neq k}{d(u_j)+d(u_k)> 1/2}}^M \biggl|\frac
{S M}{N}
\biggl(\frac{j}{M}-\frac{k}{M} \biggr) \biggr|^{2d(u_j)+2d(u_k)-2}\\
&&\qquad \leq K
\biggl(\frac{T}{N} \biggr)^{-\varepsilon}
\sum_{j,k=1}^{M} \biggl|\frac{j}{M}-\frac{k}{M} \biggr|^{\varepsilon-1}\\
&&\qquad \leq K \biggl(\frac{T}{N} \biggr)^{-\varepsilon} M^2 \int_0^1\int_0^1
|x-y|^{\varepsilon-1} \,dx \,dy \\
&&\qquad\leq K \biggl(\frac{T}{N} \biggr)^{-\varepsilon} M^2
\leq K M^2.
\end{eqnarray*}
Thus,
%
\begin{equation} \label{delta111}
\bigl|\Delta_N^{(1.1.1)} \bigr| \leq K \biggl(\frac{N}{S} \biggr)^{2} M \log M +
M^2.
\end{equation}
On the other hand, if $z= \frac{S}{N}(k-j)$ then $0<z\leq1+\delta$
for $j,k \in I_2(\delta)$. Thus, an application of Lemma \ref{z}
yields for $2>\delta>0$
\[
\bigl|\Delta_N^{(1.1.2)} \bigr| \leq K \sum_{I_2(\delta)} \biggl|1- \frac
{S}{N}(k-j) \biggr|^{2d-1},
\]
where $d:=\inf_{0\leq u \leq1}d(u)>0$. Hence, by defining
$p=k-j$ and $P=N/S$ we can write
\begin{eqnarray*}
\bigl|\Delta_N^{(1.1.2)} \bigr| &\leq& K M \sum_{p=1}^{P(1+\delta)} \biggl|1- \frac
{p}{P} \biggr|^{2d-1}\\
&\leq& K M \frac{N}{S} \int_0^{1+\delta} |1-x|^{2d-1} \,dx
\leq K M \frac{N}{S}.
\end{eqnarray*}
Note that from assumption \hyperlink{A4}{A3}, $N/S\to\infty$. Thus, by
combining the above bound and (\ref{delta111}) we conclude that
%
\begin{equation} \label{delta11}
\bigl|\Delta_N^{(1.1)} \bigr| \leq K \biggl(\frac{N}{S} \biggr)^{2} M \log M +
M^2.
\end{equation}

A similar bound can be found for
$ |\Delta_N^{(1.2)} |$ and consequently for $ |\Delta_N^{(1)} |$.
Furthermore, an analogous argument yields a similar bound for the term
$ |\Delta_N^{(2)} |$ appearing in (\ref{BN}). Now, we focus on
obtaining an
upper bound for the remaining term $R_N$ from~(\ref{RN}). By integrating
that expression
with respect to $\lambda$ we get
\begin{eqnarray*}
|R_N| & \leq & \frac{N}{T} \Biggl|\sum_{j,k=1}^{M}\sum_{t,s=0}^{N-1} h \biggl(
\frac{t}{N} \biggr) h \biggl( \frac{s}{N} \biggr) \widehat{\phi}_1(u_j,s-t)\\
&&\hspace*{57.2pt}{} \times\int_{\Pi}
\phi_2(u_k,\mu) A(u_j,x)A(u_k,-x)\\
&&\hspace*{81.8pt}{} \times  A(u_j,y) y^{-d(u_k)}
L_N(y+\mu) H(x-\mu) \\
&&\hspace*{81.8pt}\hspace*{24.4pt}{} \times e^{\imag x(t_j-t_k+t)+
\imag y(t_j-t_k+s)} \,dx \,dy \,d\mu\Biggr|,
\end{eqnarray*}
where the function $L_N(\cdot)$ is defined as
\[
L_N(x)= \cases{N, &\quad $|x|\leq1/N$,\cr
1/|x|, &\quad $1/N<|x|\leq\pi$.}
\]
Hence,
\begin{eqnarray*}
|R_N| & \leq &\frac{N}{T} \Biggl|\sum_{j,k=1}^{M}\sum_{t,s,p=0}^{N-1} h \biggl(
\frac{t}{N} \biggr)
h \biggl( \frac{s}{N} \biggr) h \biggl( \frac{p}{N} \biggr)\\
&&\hspace*{64.8pt}{}\times \widehat{\phi}_1(u_j,s-t)
\widehat{f}(u_j,u_k,t_j-t_k+t-p)\\
&&\hspace*{64.8pt}{} \times\int_{-\pi}^{\pi}\int_{-\pi}^{\pi} \phi_2(u_k,\mu)
A(u_j,y) y^{-d(u_k)}\\
&&\hspace*{64.8pt}\hspace*{49.3pt}{}\times L_N(\mu+y) e^{\imag y(t_j-t_k+s)+\imag p\mu}\,
dy\, d\mu\Biggr|\\
& \leq & K \frac{N}{T} \sum_{j,k=1}^{M}\sum_{t,s,p=0}^{N-1} |\widehat
{\phi}_1(u_j,s-t) |
|\widehat{f}(u_j,u_k,t_j-t_k+t-p) | \\
&&\hspace*{73.95pt}{} \times\biggl|\int_{-\pi}^{\pi}\int_{-\pi}^{\pi} \phi_2(u_k,\mu) A(u_j,y)
y^{-d(u_k)}\\
&&\hspace*{126.12pt}{}\times L_N(\mu+y) e^{\imag y(t_j-t_k+s)+\imag p\mu} \,dy \,d\mu\biggr|\\
& \leq & K \frac{N}{T} \sum_{j,k=1}^{M}\sum_{t,s,p=0}^{N-1} |\widehat
{\phi}_1(u_j,s-t) |
|\widehat{f}(u_j,u_k,t_j-t_k+t-p) | \\
&&\hspace*{73.95pt}{} \times\biggl|\int_{-\pi}^{\pi}\int_{-\pi}^{\pi} L_N(\mu+y)
y^{-d(u_j)-d(u_k)} \,dy\,
d\mu\biggr|
\\
& \leq & K \frac{N \log N}{T} \sum_{j,k=1}^{M}\sum_{t,s,p=0, t\neq
s}^{N-1} |s-t|^{-2d(u_j)-1}\\
&&\hspace*{116.78pt}{}\times
|S(j-k)+t-p|^{d(u_j)+d(u_k)-1}\\
& \leq & K \frac{N^2 \log N}{T} \sum_{j,k=1}^{M}\sum_{t,s,p=0, t\neq
s}^{N-1} |s-t|^{-2d(u_j)-1}
S^{d(u_j)+d(u_k)-1} \\
&&\hspace*{121pt}{}\times|j-k|^{d(u_j)+d(u_k)-1}\\
& \leq & K \frac{N^3 \log N}{T} M^2 \sum_{j, k=1}^M (S M)^{d(u_j)+d(u_k)-1}
\biggl|\frac{j}{M}-\frac{k}{M} \biggr|^{d(u_j)+d(u_k)-1} \frac{1}{M^2}.
\end{eqnarray*}
Since by assumption \hyperlink{A4}{A3}, $T/N^2\to0$, we conclude that
%
\begin{equation}\label{RN1}
|R_N| \leq K \frac{N^{3} M^2}{T^{2-d}} \int_{-\pi}^{\pi}\int_{-\pi
}^{\pi}
|x-y|^{2d-1} \,dx \,dy
\leq K N^{3} M^2 T^{d-2}.
\end{equation}
Thus, from (\ref{delta11}) and (\ref{RN1}), we conclude
\begin{eqnarray*}
&&
\frac{T B_N^{(1)}}{ [2\pi M H_{2,N}(0) ]^2}\\
&&\qquad= \frac{T}{ [2\pi M H_{2,N}(0) ]^2}\\
&&\qquad\quad{}\times\sum_{j,k=1}^{M} \int_{\Pi}
\phi_1(u_j,x)
A(u_j,x) A(u_k,-x) \phi_2(u_k,y)\\
&&\hspace*{79.7pt}{} \times  A(u_j,y) A(u_k,-y)
H_N(\lambda-x) H_N(x-\mu) H_N(\mu+y)\\
&&\hspace*{79.7pt}{} \times   H_N(-y-\lambda)
e^{\imag(x+y)(t_j-t_k)} \,dx \,dy \,d\lambda \,d\mu+ C_N,
\end{eqnarray*}
where
%
\begin{equation} \label{eq-CN}
C_N = \mathcal{O} \biggl(\frac{\log M}{S} + \frac{T}{N^2} + N T^{d-1} \biggr).
\end{equation}
Therefore, by assumption \hyperlink{A4}{A3} we conclude that $C_N = o(1)$. By following
successive decompositions as in (\ref{BN}), we replace $\phi
_2(u_k,y)$ by
$\phi_2(u_k,x)$, $A(u_k,-y)$ by $A(u_k,-x)$ and $A(u_j,y)$ by
$A(u_j,x)$, respectively.
Thus,
\begin{eqnarray*}
\frac{T B_N^{(1)}}{ [2\pi M H_{2,N}(0) ]^2}
&=& \frac{T}{ [2\pi M H_{2,N}(0) ]^2} \\
&&\times{}\sum_{j,k=1}^{M} \int_{\Pi}
\phi_1(u_j,x) A(u_j,x)
A(u_k,-x) \phi_2(u_k,x) A(u_j,x)\\
&&\hspace*{44.46pt}{} \times  A(u_k,-x) H_N(\lambda-x)H_N(x-\mu) H_N(\mu+y)\\
&&\hspace*{44.46pt}{} \times  H_N(-y-\lambda)
e^{\imag(x+y)(t_j-t_k)} \,dx \,dy \,d\lambda \,d\mu+ o(1).
\end{eqnarray*}
By integrating with respect to $\mu$ and $\lambda$, we get
\begin{eqnarray*}
&&\frac{T B_N^{(1)}}{ [2\pi M H_{2,N}(0) ]^2} \\[-1pt]
&&\qquad= \frac{T}{ [M
H_{2,N}(0) ]^2}\sum_{j,k=1}^{M}
\int_{-\pi}^{\pi}\int_{-\pi}^{\pi} \phi_1(u_j, x) A(u_j, x)\\[-1pt]
&&\qquad\quad\hspace*{72.54pt}\hspace*{50.4pt}{} \times
A(u_k,-x)\phi_2(u_k, x) A(u_j,x) A(u_k,-x) \\[-1pt]
&&\qquad\quad\hspace*{72.54pt}\hspace*{50.4pt}{} \times |H_{2,N}(x+y) |^2 e^{\imag(x+y)(t_j-t_k)} \,dx \,dy + o(1)\\[-1pt]
&&\qquad = \frac{T}{ [M H_{2,N}(0) ]^2} \sum_{j,k=1}^{M} \int_{-\pi}^{\pi}
\int_{-\pi}^{\pi} [\phi_1(u_j,x)f(u_j,x) ] [\phi
_2(u_k,x)f(u_k,x)]\\[-1pt]
&&\qquad\quad\hspace*{122.6pt}{} \times|H_{2N}(x+y) |^2 e^{\imag(x+y)[s(j-k)]} \,dx \,dy + o(1)\\[-1pt]
&&\qquad = \frac{T}{ [M H_{2,N}(0) ]^2} \sum_{j,k=1}^{M}
\int_{-\pi}^{\pi}\int_{-\pi}^{\pi} [\phi_1(u_j,x)f(u_j,x) ]
[\phi_2(u_k,x)f(u_k,x)]
\\[-1pt]
&&\qquad\quad\hspace*{122.6pt}{} \times|H_{2N}(z) |^2 e^{\imag z[s(j-k)]} \,dx \,dz+ o(1)\\
&&\qquad = \frac{T}{ [M H_{2,N}(0) ]^2}\\[-1pt]
&&\qquad\quad{}\times \sum_{j,k=1}^{M} \int_{-\pi}^{\pi}
[\phi_1(u_j,x)f(u_j,x) ] [\phi_2(u_k,x)f(u_k,x)] \\[-1pt]
&&\qquad\quad\hspace*{103.1pt}\hspace*{-51.4pt}{} \times\sum_{s,t=0}^{N-1} h^{2} \biggl( \frac{t}{N} \biggr)
h^{2} \biggl( \frac{s}{N} \biggr) \int_{-\pi}^{\pi} e^{\imag z[s(j-k)+t-s]} \,dx\,
dz+ o(1)\\[-1pt]
&&\qquad = \frac{2\pi T}{ [M H_{2,N}(0) ]^2} \sum_{s,t=0}^{N-1}\sum
_{\stackrel{j,k=1}{S(j-k)=s-t}}^{M} h^{2}
\biggl( \frac{t}{N} \biggr) h^{2} \biggl( \frac{s}{N} \biggr)\\
&&\qquad\quad\hspace*{130pt}{}\times\int_{-\pi}^{\pi} [\phi_1(u_j,x)f(u_j,x) ]\\
&&\qquad\quad\hspace*{130pt}\hspace*{28.8pt}{}\times[\phi_2(u_k,x)f(u_k,x)] \,dx + o(1).
\end{eqnarray*}
By assumption \hyperlink{A4}{A3}, for $S<N$ we can write
\begin{eqnarray*}
&&\frac{T}{ [2\pi M H_{2,N}(0) ]^2} B_N^{(1)}\\
&&\qquad= \frac{2\pi T}{ [M
H_{2,N}(0) ]^2}\\
&&\qquad\quad{}\times \sum_{t=0}^{N-1} \sum_{p=-{t/S}}^{({N-t})/{S}}
h^2 \biggl( \frac{t}{N} \biggr)
h^2 \biggl( \frac{t}{N}+\frac{p S}{N} \biggr)\\
&&\qquad\quad\hspace*{113.4pt}\hspace*{-51.8pt}{} \times\sum_{j=1}^{M-|p|}\int_{-\pi}^{\pi} [\phi_1(u_j,x)
f(u_j,x) ]\\
&&\qquad\quad\hspace*{169.2pt}\hspace*{-51.8pt}{}\times
[\phi_1(u_{j+p},x) f(u_{j+p},x) ] \,dx + o(1).
\end{eqnarray*}
Observe that by the assumptions of this proposition the products $\phi
_1(u,x)f(u,x)$
and $\phi_2(u,x) f(u,x)$ are differentiable with respect to $u$.
Furthermore, note
that by assumption \hyperlink{A4}{A3}, $\lim_{T,S\to\infty}\frac{S|p|}{T} =
0$ for any $|p|\leq
\frac{N}{S}$. Consequently,
\begin{eqnarray*}
&&\frac{S}{T}\sum_{j=1}^{M-|p|} \int_{-\pi}^{\pi} [\phi_1(u_j,x)
f(u_j,x) ]
[\phi_2(u_{j+p},x) f(u_{j+p},x) ] \,dx \\
&&\qquad \to\int_0^1\int_{-\pi}^{\pi} \phi_1(u,x) \phi_2(u,x) f(u,x)^2 \,dx
\,du,
\end{eqnarray*}
for any $|p|<\frac{N}{S}$ as $M,N,S,T\to\infty$. On the other hand,
\begin{eqnarray*}
&&\frac{2\pi T^2 N^2}{S^2 [M H_{2,N}(0)]^2}
\sum_{t=0}^{N-1} \sum_{p=-{t/S}}^{({N-t})/{S}}
h^2 \biggl( \frac{t}{N} \biggr)
h^2 \biggl( \frac{t}{N}+\frac{pS}{N} \biggr) \frac{S}{N^2}\\
&&\qquad\to
2\pi\int_0^1 \int_{-x}^{1-x} h^2(x) h^2(x+y) \,dx \,dy \biggl( \int_0^1
h^2(x) \,dx \biggr)^{-2}= 2\pi,
\end{eqnarray*}
as $M,N,S,T\to\infty$. Therefore, in this case
\[
\frac{T}{ [2\pi M H_{2,N}(0) ]^2} B_N^{(1)} \to2\pi
\int_0^1\int_{-\pi}^{\pi}\phi_1(u,x) \phi_2(u,x) f(u,x)^2 \,dx \,du,
\]
as $M,N,S,T\to\infty$. Similarly, we have that
\[
\frac{T}{ [2\pi M H_{2,N}(0) ]^2} B_N^{(2)} \to2\pi
\int_0^1\int_{-\pi}^{\pi} \phi_1(u,x) \phi_2(u,x) f(u,x)^2 \,dx \,du,
\]
as $M,N,S,T\to\infty$. Therefore, by virtue of (\ref{Tcov}) this
proposition is
proved.
\end{pf}
%
%
\begin{proposition} \label{cumulant} 
Let $\operatorname{cum}_p(\cdot)$ be the $p$th order cumulant with $p \geq
3$. Then, $T^{p/2}
\operatorname{cum}_p (J_T(\phi)) \to0$, as $ T \to\infty$.
\end{proposition}
\begin{pf} Observe that $J_T(\phi)$ can be written as
\[
J_T(\phi) = \frac{1}{2\pi M H_{2,N}(0)} Y'Q(\phi)Y,
\]
where the block-diagonal matrix $Q(\phi)$ is defined in (\ref{Q}) and
$Y \in
\mathbb{R}^{NM}$ is a Gaussian random vector defined by $Y = (Y(u_1)',
\ldots,
Y(u_M)')'$, $Y(u) = (Y_1(u)$,$\ldots, Y_N(u))$, $Y_t(u) = h
( \frac{t}{N} ) Y_{[uT] - {N/2}+t+1,T}$ with $Y_{[uT] -
{N/2}+t+1,T}$ satisfying (\ref{local-stat}). For simplicity,
denote the matrix
$Q(\phi)$ as $Q$. Since $Y$ is Gaussian,
\[
\cum_p [J_T(\phi) ] = \frac{2^{p-1}(p-1)!}{(2\pi M H_{2,N}(0))^p}
\tr(RQ)^p,
\]
where $R = \var(Y)$. Let $|A|=[\tr(AA')]^{1/2}$ be the Euclidean norm
of matrix $A$ and
let $\|A\|=\sup_{\|x\|=1} (Ax)'Ax$ be the spectral norm of $A$. Now,
since $| \tr
(QB)|\leq|Q||B|$ and $|QB| \leq\Vert Q\Vert |B|$ we get $|{\tr}(RQ)^p| \leq
\Vert RQ\Vert^{p-2}
|RQ|^2$.

On the other hand, for fixed $\lambda$, decompose the function $\phi
(\cdot,\lambda)$ as
$\phi(\cdot,\lambda)=\phi_+(\cdot,\lambda)-\phi_-(\cdot,\lambda
)$ where
$\phi_+(\cdot,\lambda) , \phi_-(\cdot,\lambda)\geq0$. Thus, we
can write
$Q=Q(\phi)=Q(\phi_+-\phi_-)=Q(\phi_+)-Q(\phi_-):=Q_+-Q_-$, say.
Now, by Lemma \ref{Lemma
110} we conclude that
\[
\|RQ\| = \|RQ_+ - RQ_-\| \leq\| RQ_+\|+\|RQ_-\|\leq K (MN^{1-2d} T^{2d-1}),
\]
and by Proposition \ref{Cov} we have that $|RQ|^2 \leq K\frac{M^2
N^2}{T}$. Thus,
\[
| {\tr}(RQ)^p |\leq K (M N^{1-2d}T^{2d-1} ) \frac{M^2 N^2}{T}.
\]
Consequently,
\[
|T^{p/2} \operatorname{cum}_p [J_T(\phi) ]|\leq K M^{1-{p/2}}
\biggl(\frac{N}{T} \biggr)^{(1-2d)(p-2)} \biggl(\frac{\sqrt{T}}{N} \biggr)^{p-2}.
\]
Since $p\geq3$ and by assumption \hyperlink{A3}{A2}, $N/T \to0$ and $\sqrt
{T}/N \to0$ as $T,N
\to\infty$, the required result is obtained.
\end{pf}

\subsection{Proof of theorems}

\mbox{}

\begin{pf*}{Proof of Theorem \ref{CON}}
To prove the consistency of the Whittle estimator, it suffices
to show that
\[
{\sup_{\theta}} |\mathcal{L}_T(\theta) - \mathcal{L}(\theta)|\to0,
\]
in probability, as $T\to\infty$, where $\mathcal{L}(\theta) := \frac{1}
{4\pi } \int_{0}^{1}\int_{-\pi}^{\pi}
[\log f_{\theta}(u,\lambda ) + \frac{f_{\theta_0} (u, \lambda)}
{f_{\theta} (u, \lambda)}]\,d\lambda\, du$. Define
$g_{\theta}(u, \lambda) = f_{\theta}(u,\lambda)^{-1}$.
By assumption \hyperlink{A1}{A1}, $g_{\theta}(u,\lambda)$ is continuous in
$\theta$, $\lambda$
and $u$. Thus, $g_{\theta}$ can be approximated by the Cesaro sum of
its Fourier series
\begin{eqnarray*}
g^{(L)}_{\theta}(u,\lambda) &=& \frac{1}{4\pi^2} \sum
_{\ell= -L}^{L}\sum_{m
= -L}^{L}
\biggl(1-\frac{|\ell|}{L} \biggr) \biggl(1-\frac{|m|}{L} \biggr)\\
&&\hspace*{70.2pt}{}\times\wihat{g}_{\theta}(\ell
,m)\exp(-i 2\pi u_j \ell- i \lambda
m),
\end{eqnarray*}
such that
$\sup_{\theta}|g_{\theta}(u,\lambda)-g^{(L)}_{\theta
}(u,\lambda)|<\varepsilon$;
see, for example, Theorem 1.5(ii) of \citet{Korn88}. Following
Theorem 3.2 of \citet{Dahl97},
we can write
\begin{eqnarray*}
&&\sup_{\theta}|\mathcal{L}_T (\theta) -\mathcal{L}(\theta)|\\
&&\qquad \leq
\mathcal{O} \biggl(\frac{1}{M} \biggr) +
\frac{\varepsilon}{4\pi}\frac{1}{M}\sum_{j=1}^{M}\int_{-\pi
}^{\pi}
[I_N(u_j,\lambda)
+f(u_j,\lambda) ] \,d\lambda\\
&&\qquad\quad{} + \frac{1}{16\pi^3} \sum_{\ell= -L}^{L}\sum_{m = -L}^{L}
\biggl(1-\frac{|\ell|}{L} \biggr) \biggl(1-\frac{|m|}{L} \biggr)\sup_{\theta}|\wihat
{g}_{\theta}(\ell,m)|\\
&&\qquad\quad\hspace*{88.4pt}{} \times\Biggl|\frac{1}{M}\sum_{j=1}^{M} \int_{-\pi}^{\pi} \exp(-i
2\pi u_j
\ell- i \lambda m)\\
&&\qquad\quad\hspace*{152.5pt}{} \times\{I_N(u_j,\lambda)-f(u_j,\lambda) \} \,d\lambda\Biggr|,
\end{eqnarray*}
where
\[
\wihat{g}_{\theta}(\ell,m) = \frac{}{} \int_{0}^{1}\int_{-\pi
}^{\pi}
g_{\theta}(u,\lambda) \exp(i 2\pi u \ell+ i \lambda m) \,du \,d\lambda.
\]
Consequently, $|\wihat{g}_{\theta}(\ell,m)|\leq
2\pi\sup_{(\theta,u,\lambda)}|g_{\theta}(u,\lambda)|$. However,
by assumption \hyperlink{A1}{A1},
$|g_{\theta}(u,\lambda)|$ is continuous in $\theta$, $u$ and
$\lambda$. Thus, since the
parameter space is compact we have that $|\wihat{g}_{\theta}(\ell
,m)|\leq K$, for some
positive constant $K$. Now, by defining for fixed $\ell, m=1,\ldots
,L$, $\phi(u,\lambda)=\cos(2\pi u
\ell)\cos(\lambda m)$ or $\phi(u,\lambda)=\sin(2\pi u \ell)\cos
(\lambda m)$ in
Proposition \ref{E} and $\phi_1(u,\lambda)=\phi_2(u,\lambda)=\cos
(2\pi u
\ell)\cos(\lambda\times\break m)$ or $\phi_1(u,\lambda)=\phi_2(u,\lambda
)=\sin(2\pi u
\ell)\cos(\lambda m)$ in Proposition \ref{Cov}, we deduce that
%
\begin{eqnarray}\label{one}
&&\frac{1}{16\pi^3} \sum_{\ell= -L}^{L}\sum_{m =
-L}^{L} \biggl(1-\frac{|\ell|}{L} \biggr) \biggl(1-\frac{|m|}{L} \biggr)
\sup_{\theta}|\wihat{g}_{\theta}(\ell,m)| \nonumber\\
&&\qquad\quad\hspace*{42.5pt}{}\times\Biggl|\frac{1}{M}\sum_{j=1}^{M} \int_{-\pi}^{\pi} \exp(-i
2\pi u_j
\ell- i \lambda m)\nonumber\\
&&\qquad\quad\hspace*{107pt}{}\times [ I_N(u_j,\lambda)-f(u_j,\lambda) ] \,d\lambda\Biggr|
\nonumber\\
&&\qquad\leq\frac{1}{16\pi^3} \sum_{\ell
=-L}^{L}\sum_{m=-L}^{L} \biggl(1-\frac{|\ell|}{L} \biggr)
\biggl(1-\frac{|m|}{L}
\biggr)\sup_{\theta}|\wihat{g}_{\theta}(\ell,m)|\nonumber\\[-8pt]\\[-8pt]
&&\qquad\quad\hspace*{75pt}{} \times\Biggl\{ \Biggl|\frac{1}{M}\sum_{j=1}^{M} \int_{-\pi}^{\pi} \cos
(2\pi u_j
\ell)\cos(\lambda m)\nonumber\\
&&\qquad\quad\hspace*{146.6pt}{}\times [ I_N(u_j,\lambda)-f(u_j,\lambda) ]
\,d\lambda\Biggr| \nonumber\\
&&\qquad\quad\hspace*{94pt}{} + \Biggl|\frac{1}{M}\sum_{j=1}^{M} \int_{-\pi}^{\pi} \sin(2\pi u_j
\ell)\cos(\lambda m)\nonumber\\
&&\qquad\quad\hspace*{159.3pt}{}\times [ I_N(u_j,\lambda)-f(u_j,\lambda) ] \,d\lambda\Biggr|
\Biggr\}\to0\nonumber
\end{eqnarray}
and
%
\begin{equation} \label{two}
\frac{1}{M}\sum_{j=1}^{M}\int_{-\pi}^{\pi}
\{I_N(u_j,\lambda)+f(u_j,\lambda) \} \,d\lambda\to2\int_0^1
\int_{-\pi}^{\pi} f(u,\lambda)\,d\lambda \,du,
\end{equation}
in probability, as $M\to\infty$. Now, from the limits (\ref{one})
and (\ref{two}), this
theorem follows.
\end{pf*}
\begin{pf*}{Proof of Theorem \ref{CLT}}
Let $\widehat{\theta}_T$ be the parameter value that
minimizes the Whittle log-likelihood function $\mathcal{L}_T(\theta)$
given by
(\ref{whittle-like}) and let $\theta_0$ be the true value of the
parameter. By the mean
value theorem, there exists a vector $\bar{\theta}_T$ satisfying
$\|\bar{\theta}_T-\theta_0\|\leq\|\widehat{\theta}_T-\theta_0\|
$, such that
%
\begin{equation} \label{taylor}
\nabla\mathcal{L}_T(\widehat{\theta}_T)-\nabla
\mathcal{L}_T(\theta_0)= [\nabla^2\mathcal{L}_T(\bar{\theta}_T) ]
(\widehat{\theta}_T-\theta_0 ).
\end{equation}
Therefore, it suffices to show that (a)
$\nabla^2\mathcal{L}_T(\theta_0)\to\Gamma(\theta_0)$,
as $T\to\infty$; (b) $\nabla^2 \mathcal{L}_T(\bar{\theta}_T)
-\nabla^2\mathcal{L}_T(\theta_0)\to0$ in probability, as $T\to
\infty$; and (c)
$\sqrt{T} \nabla\mathcal{L}_T(\theta_0)\to
N[0,\Gamma(\theta_0)]$, in distribution, as $T\to\infty$. To this
end, observe that
\begin{eqnarray*}
\nabla^2 \mathcal{L}_T(\theta) &=& \frac{1}{4\pi} \frac1M
\sum_{j=1}^M \int_{-\pi}^{\pi}
[I_N(u_j,\lambda)-f_{\theta}(u_j,\lambda)] \nabla^2
f_{\theta}(u_j,\lambda)^{-1} \\
& &{} - \nabla f_{\theta}(u_j,\lambda)
[\nabla f_{\theta}(u_j,\lambda)^{-1} ]' \,d\lambda\\
& = &\frac{1}{4\pi} \frac1M \Biggl\{
\sum_{j=1}^M \int_{-\pi}^{\pi}
\phi(u_j,\lambda) [I_N(u_j,\lambda)-f_{\theta}(u_j,\lambda)]\\
&&\hspace*{34pt}{} + \sum_{j=1}^M \int_{-\pi}^{\pi} \nabla\log f_{\theta
}(u_j,\lambda)
[\nabla\log f_{\theta}(u_j,\lambda) ]' \,d\lambda\Biggr\}\\
&=& \frac{1}{4\pi} [J_T(\phi)-J(\phi)]+\Gamma(\theta) +\mathcal
{O} \biggl(\frac1M \biggr),
\end{eqnarray*}
%
where $\phi(u,\lambda)=\nabla^2 f_{\theta}(u,\lambda)^{-1}$.
Hence, an application of
Proposition \ref{E} and Proposition \ref{Cov} yields parts (a) and
(b). On the other
hand, part (c) can be proved by means of the cumulant method. That is,
by showing that
all the cumulants of $\sqrt{T} \nabla\mathcal{L}_T(\theta_0)$
converge to zero, excepting the second order cumulant. To this end,
note that
%
\begin{eqnarray}\label{aa}
\nabla\mathcal{L}_T(\theta_0) &=& \frac{1}{4\pi} \frac1M
\sum_{j=1}^M \int_{-\pi}^{\pi}
[I_N(u_j,\lambda)-f_{\theta_0}(u_j,\lambda)] \nabla
f_{\theta_0}(u_j,\lambda)^{-1} \,d\lambda\nonumber\\
&=& \frac{1}{4\pi} J_T(\phi)- \frac{1}{4\pi} \sum_{j=1}^M \int
_{-\pi}^{\pi}
f_{\theta_0}(u_j,\lambda) \nabla
f_{\theta_0}(u_j,\lambda)^{-1} \,d\lambda\\
&=& \frac{1}{4\pi} [J_T(\phi)-J(\phi)]+\mathcal{O} \biggl(\frac{1}{M}
\biggr),\nonumber
\end{eqnarray}
where $\phi(u,\lambda)=\nabla f_{\theta_0}(u,\lambda)^{-1}$.
Hence, by Proposition \ref{E} and assumption \hyperlink{A4}{A3},
the first-order cumulant of
$\sqrt{T} \nabla\mathcal{L}_T(\theta_0)$ satisfies
%
\begin{eqnarray} \label{first-order-cum}
\sqrt{T} \esp[\nabla\mathcal{L}_T(\theta_0)]&=& \mathcal{O} \biggl(\frac
{\sqrt{T}
\log^2N}{N} \biggr) + \mathcal{O} \biggl(\frac{\sqrt{T}
}{M} \biggr) \nonumber\\[-8pt]\\[-8pt]
&\to&0,\nonumber
\end{eqnarray}
as $T\to\infty$. Furthermore, by (\ref{aa}) we have that the
second-order cumulant of
$\sqrt{T} \nabla\mathcal{L}_T(\theta_0)$ can be written as
\[
T \cov[\nabla\mathcal{L}_T(\theta_0),\nabla
\mathcal{L}_T(\theta_0)]= \frac{1}{16\pi^2}T\cov[J_T(\phi
),J_T(\phi)].
\]
Therefore, by Proposition \ref{Cov} we have that
\begin{eqnarray*}
&&\lim_{T\to\infty} T \cov[\nabla\mathcal{L}_T(\theta_0),\nabla
\mathcal{L}_T(\theta_0)] \\
&&\qquad =\frac{1}{4\pi} \int_0^1\int_{-\pi}^{\pi} \nabla f_{\theta
_0}(u,\lambda)^{-1}
[\nabla f_{\theta_0}(u,\lambda)^{-1}]'f_{\theta_0}(u,\lambda)^2
\,d\lambda \,du\\
&&\qquad =\frac{1}{4\pi} \int_0^1\int_{-\pi}^{\pi} \nabla\log f_{\theta
_0}(u,\lambda)
[\nabla\log f_{\theta_0}(u,\lambda)]' \,d\lambda \,du=\Gamma(\theta_0).
\end{eqnarray*}
Finally, for $p>2$, Proposition \ref{cumulant} gives $T^{p/2} \cum
_p[\nabla
\mathcal{L}_T(\theta_0)] \to0$, as $T\to\infty$, proving part (c).
\end{pf*}
\begin{pf*}{Proof of Theorem \ref{EFF}}
By observing that the Fisher information matrix
evaluated at the true parameter, $\Gamma_T(\theta_0)$, is given by
\[
\Gamma_T(\theta_0)=T\cov[\nabla\mathcal{L}_T(\theta_0),\nabla
\mathcal{L}_T(\theta_0)],
\]
the result is an immediate consequence of Proposition \ref{Cov}.
\end{pf*}
\begin{pf*}{Proof of Theorem \ref{average-d}}
Let $V^{(T)} = [V_{ij}^{(T)}]_{i,j=1,\ldots,p}=
\var(\wihat{\beta})$, then
\begin{eqnarray*}
\int_{0}^{1} \var[\wihat{d}(u)] \,du & = &\int_{0}^{1} \sum_{i =
1}^{p}\sum_{j = 1}^{p}
g_i(u) V_{ij}^{(T)} g_j(u) \,du \\
&=& \sum_{i = 1}^{p}\sum_{j = 1}^{p} V_{ij}^{(T)} \int_{0}^{1}
g_i(u) g_j(u) \,du\\
&=& \sum_{i = 1}^{p}\sum_{j = 1}^{p} V_{ij}^{(T)} b_{ij},
\end{eqnarray*}
where $b_{ij} = \int_{0}^{1} g_i(u) g_j(u) \,du = b_{ji}$. Therefore, by
Theorem \ref{CLT}
\[
\lim_{T\to\infty} T \int_0^1 \var[\wihat{d}(u)] \,du = \sum_{i =
1}^{p}\sum_{j =
1}^{p} \lim_{T\to\infty} \bigl[T V_{ij}^{(T)}\bigr] b_{ij}
= \sum_{i = 1}^{p}\sum_{j = 1}^{p} a_{ij} b_{ij},
\]
where $A = (a_{ij})_{i,j=1,\ldots,p}=\Gamma^{-1}$ and
\[
\Gamma_{ij} =
\frac{1}{4\pi} \int_0^1 \int_{-\pi}^{\pi} \frac{\partial
}{\partial\beta_i}
\log f(u,\lambda) \,\frac{\partial}{\partial\beta_j} \log
f(u,\lambda) \,d\lambda\,
du.
\]
But, $\log f(u,\lambda) = \log(\sigma^2)-\log(2\pi)-{d_{\beta
}(u)\log}|1-e^{\imag
\lambda}|^2$.
Thus,
\[
\frac{\partial}{\partial\beta_i} \log f(u,\lambda) = {-g_i(u)
\log}|1-e^{\imag\lambda}|^2.
\]
Hence, $\Gamma_{ij} = \int_0^1 g_i(u)g_j(u) \,du \times\frac{1}{4\pi
} \int_{-\pi}^{\pi}
(\log|1-e^{\imag\lambda}|^2 )^2 \,d\lambda= \frac{\pi^2}{6} b_{ij}$.
Therefore, $\Gamma= \frac{\pi^2}{6} B$ and $A =
\frac{6}{\pi^2} B^{-1}$. Consequently, since $A$ and $B$ are symmetric
matrices
$\lim_{T\to\infty} T \int_0^1 \var[\wihat{d}(u)] \,du = \tr(AB) =
\frac{6}{\pi^2}
\tr(I_p) = \frac{6p}{\pi^2}$. %
\end{pf*}

\section{Simulations} \label{simulations}
\begin{figure}

\includegraphics{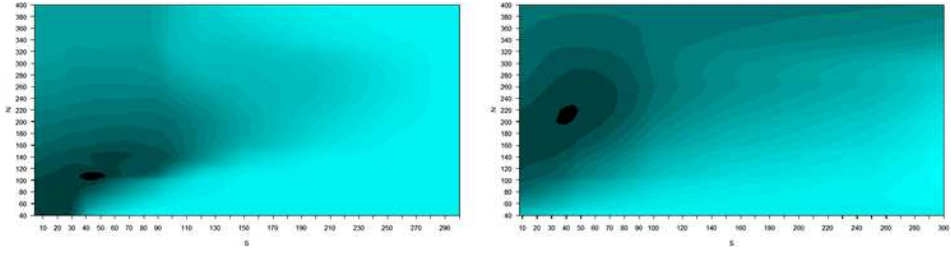}

\caption{Contour curves of the empirical MSE of Whittle
estimator. Left: sample size $T=512$. Right: sample size
$T=1024$.}\label{contorno512M2gris}
\end{figure}

In order to gain some insight into the finite sample performance of the Whittle
estimator discussed in Section \ref{methodology}, we report next a
number of Monte
Carlo experiments for the LSARFIMA model 
\[
Y_{t,T}=\sigma(t/T) (1-\vartheta B) (1-B)^{-d(t/T)} \ve_{t},
\]
for $t=1,\ldots,T$ with $d(u)=\alpha_0+\alpha_1 u$, $ \sigma
(u)=\beta_0+\beta_1 u$
and Gaussian white noise $\{\ve_{t}\}$ with unit variance. The samples
of this LSARFIMA
process are generated by means of the innovation algorithm; see, for example,
Brockwell and Davis (\citeyear{Broc91}), page 172. In this implementation, the covariances of the
process $\{Y_{t,T}\}$ is given by
\begin{eqnarray*}
\esp[Y_{s,T}Y_{t,T}]
&=&\sigma\biggl( \frac{s}{T} \biggr)\sigma\biggl( \frac{t}{T} \biggr)
\frac{\Gamma[1-d ( {s/T} )-d ( {t/T} ) ]
\Gamma[s-t+d ( {s/T} ) ]}{\Gamma[1-d ( {s/T} ) ]
\Gamma[d ( {s/T} ) ]
\Gamma[s-t+1-d ( {t/T} ) ]}
\\
&&{} \times\biggl[ 1+\vartheta^2
-\vartheta\frac{s-t-d ( {t/T} )}{s-t-1+d ( {s/T} )}
-\vartheta\frac{s-t+d ( {s/T} )}{s-t+1-d ( {t/T} )}
\biggr],
\end{eqnarray*}
for $s,t=1,\ldots,T$, $s\geq t$. Let $\theta=(\alpha_0,\alpha
_1,\beta_0,\beta_1,\vartheta)'$
be the parameter vector. The Whittle estimates in these Monte Carlo
simulations have been computed by using the cosine bell data taper
(\ref{bell}). Figure \ref{contorno512M2gris} displays the contour
curves for the empirical mean squared error (MSE) for the Whittle
estimator $\widehat{\theta}$ defined in this case as the average of
$\|\widehat{\theta} -\theta\|^2$ over 100 replications of $\widehat
{\theta}$, where $\theta$ is the true value of the parameter. These
contour curves correspond to $\theta=(0.20,0.25,0.5,0.3,0.5)$, for
sample sizes $T=512$ and $T=1024$, respectively. In these graphs, the
darkest regions represent the minimal empirical MSE while clear regions
indicate greater MSE values.
Note that for the case $T=512$, shown in the left panel, the minimal
empirical MSE region is located around $N\approx105$ and $S\approx
35$. For the sample size $T=1024$, displayed on the right panel, the
minimal empirical MSE is reached close to $N \approx200$ y $S\approx
45$. As noted in these graphs, there is a degree of flexibility for
selecting $N$ and $S$ as long they belong to the areas with minimal
empirical MSE.
Contour curves for other parameters $\theta$ such as those presented in
Tables \ref{table1} and \ref{table2} are similar to Figure \ref
{contorno512M2gris} and
produce similar empirical optimal regions for $N$ and $S$.
%
%
\begin{table}
\tabcolsep=0pt
\caption{Whittle estimation: sample size $T = 512$, block size
$N=105$ and shift
$S=35$} \label{table1}
\begin{tabular*}{\textwidth}{@{\extracolsep{\fill}}lcd{2.2}cd{2.1}d{2.1}cd{2.3}cd{2.3}d{2.3}@{}}
\hline
& \multicolumn{5}{c}{\textbf{Parameters}} & \multicolumn{5}{c@{}}{\textbf{Estimates}}
\\[-4pt]
& \multicolumn{5}{c}{\hrulefill} & \multicolumn{5}{c@{}}{\hrulefill}\\
\textbf{Case} & \multicolumn{1}{c}{$\bolds{\alpha_0}$} &
\multicolumn{1}{c}{$\bolds{\alpha_1}$} & \multicolumn{1}{c}{$\bolds{\beta_0}$}
& \multicolumn{1}{c}{$\bolds{\beta_1}$} & \multicolumn{1}{c}{$\bolds{\vartheta}$}
& \multicolumn{1}{c}{$\bolds{\widehat{\alpha}_0}$}
& \multicolumn{1}{c}{$\bolds{\widehat{\alpha}_1}$}
& \multicolumn{1}{c}{$\bolds{\widehat{\beta}_0}$}
& \multicolumn{1}{c}{$\bolds{\widehat{\beta}_1}$}
& \multicolumn{1}{c@{}}{$\bolds{\widehat{\vartheta}}$}\\
\hline
\phantom{0}1 & 0.15 & 0.20 & 0.5 & 0.3 & 0.5 & 0.130 & 0.177 & 0.497 & 0.299 &
0.473 \\
\phantom{0}2 & 0.15 & 0.20 & 0.8 & -0.2 & 0.5 & 0.124 & 0.167 & 0.795 & -0.201 &
0.463 \\
\phantom{0}3 & 0.20 & 0.25 & 0.5 & 0.3 & 0.5 & 0.161 & 0.219 & 0.497 & 0.301 &
0.455 \\
\phantom{0}4 & 0.20 & 0.25 & 0.8 & -0.2 & 0.5 & 0.163 & 0.218 & 0.797 & -0.201 &
0.453 \\
\phantom{0}5 & 0.30 & -0.20 & 0.5 & 0.3 & 0.5 & 0.291 & -0.183 & 0.498 & 0.299 &
0.506 \\
\phantom{0}6 & 0.30 & -0.20 & 0.8 & -0.2 & 0.5 & 0.287 & -0.183 & 0.797 & -0.203 &
0.501 \\
\phantom{0}7 & 0.15 & 0.20 & 0.5 & 0.3 & -0.4 & 0.138 & 0.189 & 0.496 & 0.301 &
-0.407 \\
\phantom{0}8 & 0.15 & 0.20 & 0.8 & -0.2 & -0.4 & 0.138 & 0.190 & 0.799 & -0.206 &
-0.410 \\
\phantom{0}9 & 0.20 & 0.25 & 0.5 & 0.3 & -0.4 & 0.195 & 0.228 & 0.498 & 0.299 &
-0.409 \\
10 & 0.20 & 0.25 & 0.8 & -0.2 & -0.4 & 0.193 & 0.229 & 0.795 & -0.197 &
-0.412 \\
11 & 0.30 & -0.20 & 0.5 & 0.3 & -0.4 & 0.286 & -0.197 & 0.498 & 0.298 &
-0.404 \\
12 & 0.30 & -0.20 & 0.8 & -0.2 & -0.4 & 0.279 & -0.180 & 0.796 & -0.203
& -0.404 \\
\hline
\end{tabular*}
\begin{tabular*}{\textwidth}{@{\extracolsep{\fill}}lcccccccccc@{}}
& \multicolumn{5}{c}{\textbf{Theoretical SD}} & \multicolumn{5}{c@{}}{\textbf{Estimated
SD}} \\[-4pt]
& \multicolumn{5}{c}{\hrulefill} & \multicolumn{5}{c@{}}{\hrulefill}\\
\textbf{Case} & \multicolumn{1}{c}{$\bolds{\sigma(\widehat{\alpha}_0)}$} &
\multicolumn{1}{c}{$\bolds{\sigma(\widehat{\alpha}_1)}$} & \multicolumn{1}{c}{$\bolds{\sigma(\widehat{\beta}_0)}$}
& \multicolumn{1}{c}{$\bolds{\sigma(\widehat{\beta}_1)}$}
& \multicolumn{1}{c}{$\bolds{\sigma(\widehat{\vartheta})}$}
& \multicolumn{1}{c}{$\bolds{\widehat{\sigma}(\widehat{\alpha}_0)}$}
& \multicolumn{1}{c}{$\bolds{\widehat{\sigma}(\widehat{\alpha}_1)}$}
& \multicolumn{1}{c}{$\bolds{\widehat{\sigma}(\widehat{\beta}_0)}$}
& \multicolumn{1}{c}{$\bolds{\widehat{\sigma}(\widehat{\beta}_1)}$}
& \multicolumn{1}{c@{}}{$\bolds{\widehat{\sigma}(\widehat{\vartheta})}$} \\
\hline
\phantom{0}1 & 0.115 & 0.119 & 0.035 & 0.069 & 0.109 & 0.117 & 0.146 & 0.045 &
0.089 & 0.106 \\
\phantom{0}2 & 0.115 & 0.119 & 0.047 & 0.075 & 0.109 & 0.115 & 0.146 & 0.057 &
0.100 & 0.103 \\
\phantom{0}3 & 0.115 & 0.119 & 0.035 & 0.069 & 0.109 & 0.107 & 0.132 & 0.043 &
0.091 & 0.096 \\
\phantom{0}4 & 0.115 & 0.119 & 0.047 & 0.075 & 0.109 & 0.110 & 0.131 & 0.056 &
0.098 & 0.102 \\
\phantom{0}5 & 0.115 & 0.119 & 0.035 & 0.069 & 0.109 & 0.131 & 0.140 & 0.043 &
0.090 & 0.108 \\
\phantom{0}6 & 0.115 & 0.119 & 0.047 & 0.075 & 0.109 & 0.125 & 0.140 & 0.057 &
0.099 & 0.107 \\
\phantom{0}7 & 0.074 & 0.119 & 0.035 & 0.069 & 0.051 & 0.089 & 0.155 & 0.045 &
0.091 & 0.058 \\
\phantom{0}8 & 0.074 & 0.119 & 0.047 & 0.075 & 0.051 & 0.088 & 0.150 & 0.054 &
0.096 & 0.053 \\
\phantom{0}9 & 0.074 & 0.119 & 0.035 & 0.069 & 0.051 & 0.090 & 0.142 & 0.044 &
0.091 & 0.053 \\
10 & 0.074 & 0.119 & 0.047 & 0.075 & 0.051 & 0.088 & 0.142 & 0.057 &
0.099 & 0.054 \\
11 & 0.074 & 0.119 & 0.035 & 0.069 & 0.051 & 0.089 & 0.140 & 0.046 &
0.093 & 0.055 \\
12 & 0.074 & 0.119 & 0.047 & 0.075 & 0.051 & 0.093 & 0.146 & 0.057 &
0.099 & 0.056 \\
\hline
\end{tabular*}
\end{table}
Tables \ref{table1} and \ref{table2} report the results from the
Monte Carlo simulations
for several parameter values, based on 1000 replications. These tables
show the
average of the estimates as well as their theoretical and empirical
standard deviations
(SD). The theoretical SD are based on Theorem \ref{CLT} with matrix
$\Gamma_{\theta}$ given by
\[
\Gamma_{\theta} = \pmatrix{
\Gamma_{\alpha} & 0& \gamma_{\alpha\vartheta}   \cr
0  & \Gamma_{\beta} & 0 \cr
\gamma_{\alpha\vartheta}' & 0 & \gamma_{\vartheta} },
\]
where
$\gamma_{\alpha\vartheta} = [\frac{\log(1-\vartheta)}{\vartheta
}, \frac{\log(1-\vartheta)}{2\vartheta} ]', \gamma_{\vartheta} =
\frac{1}{1-\vartheta^2}$,
and the matrices $\Gamma_{\alpha}$ and $\Gamma_{\beta}$ are given in
Example \ref{examp-poly}. The bandwidth parameters $N$ and $S$ for
each table are based on values found in Figure \ref{contorno512M2gris}
for $\theta=(0.20,0.25,0.5,0.3,0.5)$. As
mentioned above, these values are very similar for the other parameters
reported in
Tables \ref{table1} and \ref{table2}. Observe from these tables that
the estimated
parameters are close to their true values. Besides, the empirical
standard deviations
are close to their theoretical counterparts.
%
\begin{table}
\tabcolsep=0pt
\caption{Whittle estimation: sample size $T = 1024$, block size
$N=200$ and shift
$S=45$}\label{table2}
\begin{tabular*}{\textwidth}{@{\extracolsep{\fill}}lcd{2.2}cd{2.1}d{2.1}cd{2.3}cd{2.3}d{2.3}@{}}
\hline
& \multicolumn{5}{c}{\textbf{Parameters}} & \multicolumn{5}{c@{}}{\textbf{Estimates}}
\\[-4pt]
& \multicolumn{5}{c}{\hrulefill} & \multicolumn{5}{c@{}}{\hrulefill}\\
\textbf{Case} & \multicolumn{1}{c}{$\bolds{\alpha_0}$} &
\multicolumn{1}{c}{$\bolds{\alpha_1}$} & \multicolumn{1}{c}{$\bolds{\beta_0}$}
& \multicolumn{1}{c}{$\bolds{\beta_1}$} & \multicolumn{1}{c}{$\bolds{\vartheta}$}
& \multicolumn{1}{c}{$\bolds{\widehat{\alpha}_0}$}
& \multicolumn{1}{c}{$\bolds{\widehat{\alpha}_1}$}
& \multicolumn{1}{c}{$\bolds{\widehat{\beta}_0}$}
& \multicolumn{1}{c}{$\bolds{\widehat{\beta}_1}$}
& \multicolumn{1}{c@{}}{$\bolds{\widehat{\vartheta}}$}\\
\hline
\phantom{0}1 & 0.15 & 0.20 & 0.5 & 0.3 & 0.5 & 0.127 & 0.193 & 0.498 & 0.301 &
0.475 \\
\phantom{0}2 & 0.15 & 0.20 & 0.8 & -0.2 & 0.5 & 0.131 & 0.195 & 0.796 & -0.198 &
0.479 \\
\phantom{0}3 & 0.20 & 0.25 & 0.5 & 0.3 & 0.5 & 0.179 & 0.239 & 0.497 & 0.304 &
0.473 \\
\phantom{0}4 & 0.20 & 0.25 & 0.8 & -0.2 & 0.5 & 0.176 & 0.241 & 0.798 & -0.199 &
0.475 \\
\phantom{0}5 & 0.30 & -0.20 & 0.5 & 0.3 & 0.5 & 0.286 & -0.189 & 0.500 & 0.298 &
0.493 \\
\phantom{0}6 & 0.30 & -0.20 & 0.8 & -0.2 & 0.5 & 0.286 & -0.197 & 0.799 & -0.202 &
0.488 \\
\phantom{0}7 & 0.15 & 0.20 & 0.5 & 0.3 & -0.4 & 0.143 & 0.198 & 0.498 & 0.302 &
-0.404 \\
\phantom{0}8 & 0.15 & 0.20 & 0.8 & -0.2 & -0.4 & 0.144 & 0.197 & 0.797 & -0.197 &
-0.404 \\
\phantom{0}9 & 0.20 & 0.25 & 0.5 & 0.3 & -0.4 & 0.195 & 0.245 & 0.500 & 0.300 &
-0.405 \\
10 & 0.20 & 0.25 & 0.8 & -0.2 & -0.4 & 0.197 & 0.243 & 0.797 & -0.199 &
-0.405 \\
11 & 0.30 & -0.20 & 0.5 & 0.3 & -0.4 & 0.293 & -0.197 & 0.500 & 0.299 &
-0.402 \\
12 & 0.30 & -0.20 & 0.8 & -0.2 & -0.4 & 0.293 & -0.200 & 0.797 & -0.199
& -0.403 \\
\hline
\end{tabular*}
\begin{tabular*}{\textwidth}{@{\extracolsep{\fill}}lcccccccccc@{}}
& \multicolumn{5}{c}{\textbf{Theoretical SD}} & \multicolumn{5}{c@{}}{\textbf{Estimated
SD}} \\[-4pt]
& \multicolumn{5}{c}{\hrulefill} & \multicolumn{5}{c@{}}{\hrulefill}\\
\textbf{Case} & \multicolumn{1}{c}{$\bolds{\sigma(\widehat{\alpha}_0)}$} &
\multicolumn{1}{c}{$\bolds{\sigma(\widehat{\alpha}_1)}$} & \multicolumn{1}{c}{$\bolds{\sigma(\widehat{\beta}_0)}$}
& \multicolumn{1}{c}{$\bolds{\sigma(\widehat{\beta}_1)}$}
& \multicolumn{1}{c}{$\bolds{\sigma(\widehat{\vartheta})}$}
& \multicolumn{1}{c}{$\bolds{\widehat{\sigma}(\widehat{\alpha}_0)}$}
& \multicolumn{1}{c}{$\bolds{\widehat{\sigma}(\widehat{\alpha}_1)}$}
& \multicolumn{1}{c}{$\bolds{\widehat{\sigma}(\widehat{\beta}_0)}$}
& \multicolumn{1}{c}{$\bolds{\widehat{\sigma}(\widehat{\beta}_1)}$}
& \multicolumn{1}{c@{}}{$\bolds{\widehat{\sigma}(\widehat{\vartheta})}$} \\
\hline
\phantom{0}1 & 0.081 & 0.084 & 0.025 & 0.049 & 0.077 & 0.089 & 0.106 & 0.032 &
0.064 & 0.081 \\
\phantom{0}2 & 0.081 & 0.084 & 0.033 & 0.053 & 0.077 & 0.097 & 0.106 & 0.040 &
0.069 & 0.085 \\
\phantom{0}3 & 0.081 & 0.084 & 0.025 & 0.049 & 0.077 & 0.093 & 0.097 & 0.031 &
0.062 & 0.078 \\
\phantom{0}4 & 0.081 & 0.084 & 0.033 & 0.053 & 0.077 & 0.090 & 0.095 & 0.038 &
0.067 & 0.073 \\
\phantom{0}5 & 0.081 & 0.084 & 0.025 & 0.049 & 0.077 & 0.107 & 0.103 & 0.030 &
0.061 & 0.091 \\
\phantom{0}6 & 0.081 & 0.084 & 0.033 & 0.053 & 0.077 & 0.101 & 0.104 & 0.040 &
0.068 & 0.079 \\
\phantom{0}7 & 0.052 & 0.084 & 0.025 & 0.049 & 0.036 & 0.066 & 0.110 & 0.031 &
0.061 & 0.039 \\
\phantom{0}8 & 0.052 & 0.084 & 0.033 & 0.053 & 0.036 & 0.065 & 0.113 & 0.039 &
0.066 & 0.040 \\
\phantom{0}9 & 0.052 & 0.084 & 0.025 & 0.049 & 0.036 & 0.066 & 0.100 & 0.030 &
0.060 & 0.040 \\
10 & 0.052 & 0.084 & 0.033 & 0.053 & 0.036 & 0.058 & 0.087 & 0.040 &
0.068 & 0.037 \\
11 & 0.052 & 0.084 & 0.025 & 0.049 & 0.036 & 0.066 & 0.103 & 0.029 &
0.060 & 0.039 \\
12 & 0.052 & 0.084 & 0.033 & 0.053 & 0.036 & 0.064 & 0.101 & 0.039 &
0.068 & 0.039 \\
\hline
\end{tabular*}
\end{table}
These simulations suggest that the finite sample performance of the
proposed estimators
seems to be very good in terms of bias and standard deviations. This,
despite the fact
that in many of these simulations we have tested the method with large
values of the
long-memory parameter, that is, close to $ \frac12$. In Table \ref{table1}, for example,
for the combination $\alpha_0=0.20$, $\alpha_1=0.25$, the maximum
value of $d(u)$ is
$0.45$. Additional Monte Carlo experiments with other model
specifications are reported
in \citet{Palm10a}. Those simulations explore the empirical
optimal selection of $N$ and
$S$ and the finite sample performance of the Whittle estimators. Note,
however, that
further research is needed to establish optimal selection of $N$ and
$S$ from a
theoretical perspective. A comparison of the performances of the
Whittle method with a
kernel maximum likelihood estimation approach proposed by \citet
{Bera09} and two data
illustrations are also discussed in that paper.

\section{Final remarks} \label{conclusion}
A class of locally stationary long-memory processes has been addressed
in this paper, which is capable of modeling nonstationary time series
data exhibiting
time-varying long-range dependence. A computationally efficient Whittle
estimation
method has been proposed and it has been shown that these estimators
possess very
desirable asymptotic properties such as consistency, normality and
efficiency. Moreover,
several Monte Carlo simulations indicate that the estimates perform
well even for
relatively small sample sizes.

\begin{appendix}
\section*{Appendix}
This appendix contains nine auxiliary lemmas used to prove the theorems
stated in
Section \ref{methodology} and the propositions stated in Section \ref{proofs}. Proof
of these results
are provided in \citet{Palm10a}.
%
\begin{lemma} \label{Lemma A.1}
Let $f(u,\lambda)$ be a time-varying spectral density satisfying
assumption \textup{\hyperlink{A1}{A1}}
and let $\phi\dvtx[0,1]\times[-\pi,\pi]\to\mathbb{R}$ be a function
such that
$\phi(u,\lambda)$ is continuously differentiable in $\lambda$.
Consider the function
defined by
\[
g(u,\lambda) = \int_{-\pi}^{\pi} \phi(u, \lambda+ \omega
)f(u,\omega) \,d\omega,
\]
and its Fourier coefficients
$\widehat{g}(u,k) = \int_{-\pi}^{\pi} g(u,\lambda) e^{-\imag k
\lambda} \,d\lambda$.
Under assumption \textup{\hyperlink{A1}{A1}}, for every $u\in[0,1]$ we have that
$\lim_{n\to\infty} \sum_{k = -n}^{n} \widehat{g}(u,k) = 2\pi g(u,0)$.
\end{lemma}
\begin{lemma} \label{A.5}
Consider the function $\phi\dvtx [0,1]\times[-\pi,\pi]\to\mathbb{C}$,
such that $\partial
\phi(u$,\break $\gamma) /\partial u $ exists and $|\partial\phi(u,\gamma)
/\partial u| \leq K
|\gamma|^{-2d(u)}$, where $0\leq d(u)\leq d$ for all $u\in[0,1]$.
Then, for any $0\leq
t\leq N$ we have that
\[
H_N \biggl[\phi\biggl( \frac{\cdot}{T},\gamma\biggr)
h \biggl( \frac{\cdot}{N} \biggr),\lambda\biggr] =
\phi\biggl( \frac{t}{T},\gamma\biggr) H_N(\lambda) +
\mathcal{O} \biggl[ \frac{N}{T} |\gamma|^{-2d} L_N(\lambda) \biggr].
\]
\end{lemma}
\begin{lemma} \label{XXXX}
Consider $d_1,d_2\in[0,1/2)$ and for any $\ell\in\mathbb{Z}$ define
the integral
$I(\ell) = \int_{1}^{\infty} [(x-1)^{-2d_1} - x^{-2d_1}]
|\ell+x|^{d_1+d_2-1} \,dx$. Then $I(\ell) = \mathcal{O}(|\ell|^{d_1+d_2-1})$.
\end{lemma}
\begin{lemma}
\label{Lemma 108} Let $\phi(u, \lambda)$ be a positive function,
symmetric in $\lambda$,
such that $\phi(u, \lambda) \geq C |\lambda|^{2d(u)}$, for $\lambda
\in[-\pi, \pi]$,
where $d(u)$ is a positive bounded function for $u\in[0,1]$ and $C>0$.
Let $Q(u)$ for
$u\in[0,1]$ be the matrix defined in (\ref{Q}). Then there exists $K>
0$ such that
$X'Q(u)^{-1} X \leq K X'X N^{2d(u)}$, for all vector $X
\in\mathbb{R}^{N}$.
\end{lemma}
\begin{lemma}\label{Lemma 106}
Let $\phi(u, \lambda)$ be a positive function, symmetric in $\lambda
$, such that
$\phi(u, \lambda) \geq C |\lambda|^{2d(u)}$, for $\lambda\in[-\pi
, \pi]$, where
$d(u)$ is a positive bounded function for $u\in[0,1]$ and $C>0$. Let
$Q(u)$ for
$u\in[0,1]$ and $Q(\phi)$ be the matrices defined in (\ref{Q}). Then
there exists $K>
0$ such that
\[
|X'[Q(\phi)^{-1}-Q(\varphi)]X | \leq K X'X N^{2d + {1/2}},
\]
where $\varphi(u,\cdot) = \phi(u, \cdot)^{-1}/4 \pi^2$, $d = \sup
d(u)<\infty$ and $X
\in\mathbb{R}^{NM}$.
\end{lemma}
\begin{lemma} \label{Lemma 110}
Let $\phi(u, \lambda)$ be a positive function, symmetric in $\lambda
$, such that
$\phi(u, \lambda) \geq C |\lambda|^{2d(u)}$, for $\lambda\in[-\pi
, \pi]$, where
$d(u)$ is a positive bounded function for $u\in[0,1]$ and $C>0$. Let
$Q(\phi)$ be the
block-diagonal matrix defined in (\ref{Q}). Then there exists $K> 0$
such that
\[
\sup_X \biggl|\frac{X'RX}{X'Q(\phi)^{-1}X} \biggr| \leq K M N^{1-2d}
T^{2d-1},
\]
where $d = \sup d(u)<\frac{1}{2}$ and $X \in\mathbb{R}^{NM}$.
\end{lemma}
\begin{lemma}
\label{Lemma 101} Let $f(\lambda)$ and $\phi(\lambda)$ be two
real-valued functions
defined over $\lambda\in[-\pi,\pi]$ with Fourier coefficients
$\widehat{f}(k)$ and
$\widehat{\phi}(k)$, respectively, satisfying $|\widehat
{f}(k)\widehat{\phi}(k)|\leq
K/k^2$, for some positive constant $K$ and $|k|>0$. Let $C(N)$ be given
by $C(N) =
\sum_{t=0}^{N-1} h^2 ( \frac{t}{N} ) \sum_{k=N-t}^{N-1}
\widehat{f} (k) \widehat{\phi}(k)$ with bounded data taper,
$|h(u)|<K$, for all
$u\in[0,1]$. Then there exits a positive constant $K$ such that
$|C(N)|\leq K \log^2
N$.
\end{lemma}
\begin{lemma} \label{Lemma 102}
Define $D(N,T) = \frac{1}{N} \sum_{t=0}^{N-1} \sum_{k=N-t+1}^{N-1}
\frac{\varphi(k)}{k^2-d^2}
( \frac{t-N/2}{T} )$
with function $|\varphi(k) |< C \log N$ for all $0 \leq k \leq N$,
$N>1$, where $C$ is a
positive constant. Then there exists a constant $K>0$ such that
$|D(N,T)| \leq K
\frac{\log^2 N}{T}$.
\end{lemma}
\begin{lemma} \label{z}
Let $z\in[0,1+\delta]$ with $2>\delta>0$ and $2\beta>2\alpha>0$.
Then, the positive
double integral $ I(z)=\int_0^1 |z-x|^{\alpha-1} \int_1^\infty
(y-x)^{-\beta} (y-z)^{\alpha-1} \,dy \,dx$,
satisfies $I(z)\leq K |1-z|^{2\alpha-\beta}$.
\end{lemma}
\end{appendix}

\section*{Acknowledgments}
We are deeply thankful to the Associate Editor and two anonymous
referees for their careful reading of the manuscript and for their
constructive comments which led to substantial improvements.

\printaddresses

\end{document}